\documentclass[twoside,12pt,reqno]{amsart}
\usepackage{amsmath}
\usepackage{amscd}
\usepackage{latexsym}
\theoremstyle{plain}
\newtheorem{thm}{Theorem}[section]
\newtheorem{cor}[thm]{Corollary}
\newtheorem{lem}[thm]{Lemma}
\newtheorem{prop}[thm]{Proposition}
\newtheorem{Claim}[thm]{Claim}

\theoremstyle{definition}
\newtheorem{defn}[thm]{Definition}

\newtheorem{Remark}[thm]{Remark}

\newcommand{\Hdg}[2]{\mathcal{H}^{#1,#2}}

\newcommand{\ZZ}{\mathbb{Z}}
\newcommand{\CC}{\mathbb{C}}

\newcommand{\FF}{\mathbb{F}}

\newcommand{\PP}{\mathbb{P}}

\newcommand{\SP}{\text{Spec }}

\newcommand{\HI}[2]{\text{Hilb}_{#1}(#2)}
\newcommand{\Hi}[1]{\text{Hilb}_{#1}}
\newcommand{\GR}[3]{\text{Grass}_{#1}(#2,#3)}

\def\Sym{{\rm Sym}}
\def\pic{{\rm Pic}}

\def\P{{\mathbb P}}

\def\ps{\vspace{4pt}}

\newcommand{\lt}{\left}
\newcommand{\rt}{\right}

\newcommand{\OO}{\mathcal O}

\newcommand{\rarrowtip}{{\scriptscriptstyle -\mkern -10mu \mbox{$\scriptscriptstyle\succ$}}}
\newcommand{\dottedrightarrow}{\mathbin{\raisebox{0.30ex}{\mbox{\boldmath$\scriptscriptstyle -\mkern -1mu -\mkern -1mu \rarrowtip$}}}}

\newlength{\extraitemseplength}
\newcommand{\extraitemsep}{\rule[-\extraitemseplength]{0cm}{\extraitemseplength}}

\renewcommand{\xrightarrow}[1]{\mathbin{\stackrel{#1}{\longrightarrow}}}

\begin{document}

\title{Curves of Small Degree on Cubic Threefolds}
\author[Harris]{Joe Harris} 
\address{Department of Mathematics \\ Harvard University \\ Cambridge MA 02138}
\email{harris@math.harvard.edu}
\author[Roth]{Mike Roth} 
\address{Department of Mathematics \\ University of Michigan \\ Ann Arbor, MI 48109}
\email{mikeroth@umich.edu}
\author[Starr]{Jason Starr}\thanks{The third author was partially
supported by an NSF 
Graduate Research Fellowship and Sloan Dissertation Fellowship}
\address{Department of Mathematics \\ Massachusetts Institute of Technology \\ Cambridge MA 02139}
\email{jstarr@math.mit.edu}
\date{\today}

\maketitle

\tableofcontents

\begin{abstract}
In this article we consider the spaces $\Hdg{d}{g}(X)$ parametrizing
smooth curves of degree $d$ and genus $g$ on a smooth cubic threefold $X
\subset \PP^4$.  For $1\leq d\leq 5$, we show that each variety
$\Hdg{d}{g}(X)$ is irreducible of dimension $2d$.  
\end{abstract}

\section{Introduction} \label{sec-intro}

Suppose that $X\subset \PP^4$ is a smooth cubic hypersurface in
complex projective 4-space. In this article we consider the space
$\Hdg{d}{g}(X)$ parametrizing smooth curves of degree $d$ and genus
$g$ on a smooth cubic threefold $X\subset\PP^4$.  For $1\leq d\leq 5$
we show that each variety $\Hdg{d}{g}(X)$ is irreducible of dimension
$2d$.  

\

For the \emph{Fano scheme} of lines $F=\Hdg{1}{0}(X)$, this is a
classical result, c.f. ~\cite{CG}.  We bootstrap from this case by
residuation: in each case we show that for a general point
$[C]\in\Hdg{d}{g}(X)$ there is a surface $\Sigma\subset \PP^4$ which
contains $C$ and such that every irreducible component of the residual 
to $C$ in $\Sigma\cap X$ has degree $e<d$.  In this way we inductively 
prove that for $1\leq d\leq 5$ the space $\Hdg{d}{g}(X)$ is
irreducible, and in several cases we also show smoothness.  In a
forthcoming paper ~\cite{HRS2}, we use similar methods to describe the 
\emph{Abel-Jacobi maps} $u_{d,g}:\Hdg{d}{g}(X)\rightarrow J(X)$ for
$1\leq d\leq 5$.  

\subsection{Notation} \label{sec-not}

All schemes in this paper will be schemes over $\CC$.  All absolute
products will be understood to be fiber products over
$\text{Spec}(\CC)$.

\ps

For a projective variety $X$ and a numerical polynomial $P(t)$,
$\HI{P(t)}{X}$ denotes the corresponding Hilbert scheme.  For integers
$d,g$, $\Hdg{d}{g}(X)\subset \HI{dt+1-g}{X}$ denotes the open
subscheme parametrizing smooth, connected curves of degree $d$ and
genus $g$.

\section{Preliminaries}

In this section we gather some preliminary facts about deformation
theory, residuation, and Abel-Jacobi maps.

\ps

\subsection{Deformation Theory}~\label{sec-def}

All of the irreducibility arguments in this paper follow the same
pattern, and the linchpin of these arguments is the infinitesimal
analysis of the Hilbert scheme in~\cite[section I.2]{K}, in
particular~\cite[theorem I.2.15]{K}.  The part of this theorem which 
we shall use most often is the following:

\ps

\begin{prop}~\label{prop-def1}  
Let $Y$ be a smooth complex variety
with canonical 
divisor class $K_Y$ and let $C\subset Y$ be a connected, local complete
intersection curve with normal bundle $N_{C/Y} =
\mathcal{I}_C/\mathcal{I}_C^2$ and with arithmetic genus $p_a$.   
Every irreducible component of the Hilbert scheme 
at $[C]$ has dimension at least
\begin{equation}
\chi(N_{C/Y}) = h^0(N_{C/Y})-h^1(N_{C/Y}) = -K_Y.[C] +
(1-p_a)(\dim Y - 3). 
\end{equation}
The Zariski tangent space has dimension $h^0(N_{C/Y})$, therefore the
Hilbert scheme is smooth at $[C]$ if $h^1(N_{C/Y})=0$.
\end{prop}

Although this is technically inaccurate, we will say that the curve
$C\subset Y$ is \emph{unobstructed} if $h^1(N_{C/Y})=0$.  

\ps

Another condition closely related to smoothness of the Hilbert scheme
at $[C]$ is the question of whether deformations of $C$ \emph{smooth}
the singularities of $C$, i.e., whether or not $C$ is in the closure of
the open set parametrizing smooth curves.  Suppose that $C$ is a
\emph{nodal} curve, i.e., every singular point is formally isomorphic
to the formal neighborhood of $0\in\SP\CC[x,y]/xy$.  Then~\cite[lemma
9.2.2]{CK}  the deformation space of the nodes
$p_1,\dots,p_{\delta}$ is canonically identified with 
\begin{equation}
H^0(C,\text{Ext}^1_{\OO_C}(\Omega_C,\OO_C))~=~\oplus_{i~=~1}^\delta
T_i'\otimes T_i''
\end{equation}
where $T_i'$, $T_i''$ are the tangent spaces of the two branches of $C$
at $p_i$.  In the case that $C$ is unobstructed we have a short exact
sequence:
\begin{equation}\label{eqn-def0}\begin{CD}
H^0(C,N_{C/X}) @>>> H^0(C,\text{Ext}^1_{\OO_C}(\Omega_C,\OO_C)) @>>>
H^1(C,T_Y|_C) @>>> 0
\end{CD}\end{equation}
This calculation leads to the following:

\ps

\begin{lem}~\label{lem-def1}
When $h^1(C,T_Y|_C)=0$ the morphism from the formal neighborhood of $[C]$
in the Hilbert scheme to the deformation space of the
nodes is smooth at $[C]$, thus deformations of $C$ smooth the nodes.
\end{lem}

\ps

A different approach to smoothing nodes is as
follows (at some level it is equivalent to the last
paragraph).  Suppose $Y$ is a smooth variety and $Z\subset Y$ is
a simple normal crossings subscheme with no triple points, in
particular each irreducible component of $Z$ is smooth.  Let $Z_i$
be an irreducible component of $Z$ and let $D_1,\dots,D_r$ be the
connected components of $\text{sing}(Z)\cap Z_i$.  For each
$j~=~1,\dots,r$, let $Z'_j$ be the second irreducible component of $Z$
which contains $D_i$ (if an irreducible component intersects itself,
make an \'etale base change such that the preimage of $Z_i$ decomposes 
into a union of irreducible components in a neighborhood of the
preimage of $D_i$).  Consider the diagram of sheaves:
\begin{equation}\begin{CD}
0 @>>>\lt(I_{Z/Y}\rt)/\lt(I_{Z/Y}I_{Z_i/Y}\rt) @>>>
\lt(I_{Z_i/Y}\rt)/\lt(I^2_{Z_i/Y}\rt) @>>>
\lt(I_{Z_i/Y}\rt)/\lt(I_{Z/Y}\rt) @>>> 0
\end{CD}\end{equation}
where $I_{A/B}$ is the ideal sheaf of $A$ in $B$.
By passing to formal neighborhoods and using the canonical form for a
simple normal crossings variety, one sees that this is a short exact
sequence.  Moreover one can identify the last term with
$\oplus_{j~=~1}^r I_{D_j/Z'_j}/I^2_{D_j/Z'_j}$.  Dualizing this short
exact sequence leads to the short exact sequence:
\begin{equation}\label{eqn-def}\begin{CD}
0 @>>> N_{Z_i/Y} @>>> N_{Z/Y}|_{Z_i} @>>> \oplus_{j~=~1}^r
N_{D_j/Z_i}\otimes N_{D_j/Z'_j} @>>> 0
\end{CD}\end{equation}
Now suppose that $Z$ is a curve with two irreducible
components $Z_1$ and $Z_2$ intersecting at a node $p$.  We have an
obvious short exact sequence of sheaves:
\begin{equation}\begin{CD}
0 @>>> N_{Z/Y}|_{Z_1}(-p) @>>> N_{Z/Y} @>>> N_{Z/Y}|_{Z_2} @>>> 0,
\end{CD}\end{equation}
and the map in equation~\ref{eqn-def0} is simply the composite map
\begin{equation}\begin{CD}
H^0(Z,N_{Z/Y}) @>>> H^0(Z_2,N_{Z/Y}|_{Z_2}) @>>> T_{Z_1,p}\otimes
T_{Z_2,p}
\end{CD}\end{equation}
where the second map comes from equation~\ref{eqn-def}.  Again using
equation~\ref{eqn-def} and combining this with the long exact sequence
in cohomology associated to a short exact sequence of sheaves, we
conclude the following
\begin{lem}\label{lem-defo}
Suppose that $Z\subset X$ is a nodal curve and $Z_1,Z_2$ are two
closed nodal subcurves of $Z$ which intersect transversally in a
single point $p\in Z_1\cap Z_2$.  
Then $Z$ is unobstructed and the node of $Z$ smooths when
$H^1(Z_1,N_{Z_1/Y}(-p))~=~H^1(Z_2,N_{Z_2/Y})~=~0$.  
\end{lem}

\ps

Let us return now to the strategy of proving that $\Hdg{d}{g}(X)$ is
irreducible.  The first case will be showing that \emph{Fano scheme}
of lines $F:=\Hdg{1}{0}(X)$ is irreducible, in fact a smooth,
projective surface.  The analysis of this case is classical.  For each 
$1<d \leq 5$, we define an incidence correspondence
\begin{equation}
f_{d,g}:I^{d,g}~\rightarrow~\Hdg{d}{g}(X)
\end{equation}
parametrizing curves $C\subset X$ along with some extra data and 
such that $f_{d,g}$ is dominant of constant fiber dimension.  The
extra data will  
allow us to associate a surface $S\subset \PP^4$ which contains $C$
and such that the residual of $C$ in $S\cap X$ is made up of curves of 
strictly smaller degree.  We stratify $I_{d,g}$ according to the
behavior of the residual curve.  By studying the residual curves in
each case, we prove that there is a unique irreducible component of
$I_{d,g}$ whose image in $\Hdg{d}{g}(X)$ has dimension $\geq 2d$, and
that this image has dimension precisely $2d$.  Then it follows that
$\Hdg{d}{g}(X)$ is irreducible of dimension $2d$.

\subsection{Residuation} \label{sec-res}
In this section we review a few basic facts about residuation of
subschemes in a Gorenstein scheme.  

\ps

\begin{defn}  Suppose that $D$ is a Gorenstein scheme and
$D_1\subset D$ is a closed subscheme of codimension 0.  Let
$\mathcal{I}$ denote the ideal sheaf of $D_1$ in $D$.  Define
\begin{equation}
\mathcal{J} = (0:_{\OO_D} \mathcal{I}) =
\underline{\text{Hom}}_{\OO_D}( \OO_{D_1}, \OO_D).
\end{equation}  
Denote by
$D_2\subset D$ the closed subscheme associated to the ideal sheaf
$\mathcal{J}$.  We define $D_2\subset D$ to be the \emph{residual
subscheme} to $D_1\subset D$.  
\end{defn}

\ps

\begin{thm}~\cite[theorem 21.23]{E}~\label{thm-link}  Let $D$ be a
Gorenstein scheme, $D_1\subset D$ a codimension 0 closed subscheme.
Let $D_2\subset D$ be the residual subscheme to $D_1\subset D$.
\begin{enumerate}
\item The codimension of $D_2\subset D$ is zero and $D_2$ has no
embedded components.  If $D_1$ has no embedded components, then
$D_1\subset D$ is the residual subscheme to $D_2\subset D$.
\item If $D_1$ is Cohen-Macaulay, then $D_2$ is Cohen-Macaulay.
\item If $D_1$ is Cohen-Macaulay, then $\mathcal{J}\otimes\omega_D$ is
a canonical sheaf for $D_1$.  In particular, $D_1$ is Gorenstein iff
$\mathcal{J}$ is locally principal.
\end{enumerate}
\end{thm}

We will often be concerned with flat families of 1-cycles.  The
question arises when flatness of $D$ and $D_1$ over $B$ implies that
$D_2$ is also flat over $B$.  The following lemma addresses this issue 
and also establishes a \emph{base-change result} for residual subschemes.

\ps

\begin{lem}~\label{lem-flatlink}  Let $R$ be a local Noetherian ring.  Let
$A$ be a local Noetherian A-algebra (i.e., $R~\rightarrow~ A$ is a
local homomorphism) such that $A$ is Gorenstein and flat over $R$.
Let $I\subset A$ be a codimension zero ideal such that $A/I$ is
Cohen-Macaulay.  Define $J=(0:_A I)$.
\begin{enumerate}
\item For any regular sequence $(r_1,\dots,r_n)$ in $R$, we have 
\begin{equation}\label{eqn-10}
 J/(r_1,\dots,r_n)J = \lt(0:_{A/(r_1,\dots,r_n)A} I/(r_1,\dots,r_n)I 
\rt).\end{equation}
\item If $R$ is regular, then $A/I$ and $A/J$ are flat over $R$.
\end{enumerate}
\end{lem}

\begin{proof}
First we prove (1).  Since $I\subset A$ has codimension zero and $A/I$ 
is Cohen-Macaulay, $A/I$ is a maximal Cohen-Macaulay module.  Since
$A$ is flat over $R$, $(r_1,\dots,r_n)$ is a regular sequence for
$A$.  Using~\cite[proposition 18.13]{E}, the result follows by induction on 
$n$ with~\cite[proposition 21.12(b)]{E} as the induction step.

\ps

Now we prove (2).  By (2) of theorem~\ref{thm-link}, we know that
$A/J$ is Cohen-Macaulay.  By~\cite[theorem 18.16]{E}, $A/J$ is flat
over $R$ iff
\begin{equation}\label{eqn-11}
\dim (A/J) = \dim (R) + \dim (A/(J+m_R A)).
\end{equation}
We always have the inequality
\begin{equation}\label{eqn-12}
\dim (A/J) \leq \dim (R) +
\dim (A/(J+m_R A)).\end{equation} 
We also have the inequality
\begin{equation}\label{eqn-13}
\dim (R) + \dim (A/(J+m_R A)) \leq \dim (R) +
\dim (A/m_R A).\end{equation}
Now $A$ is flat over $R$, so we have 
\begin{equation}\label{eqn-14}
\dim (R)+\dim (A/m_R A) = \dim (A).\end{equation}
Finally, since $J\subset A$ has codimension zero,
$\dim (A)=\dim (A/J)$.  Putting the inequalities together,
we have
\begin{equation}\label{eqn-15}
\dim (A/J)\leq \dim (R) + \dim (A/(J+m_R A)) \leq
\dim (A/J).\end{equation}
Thus $A/J$ is flat over $R$.  By the same argument $A/I$ is also flat
over $R$.
\end{proof}

\ps

\begin{cor}[Reformulation]~\label{cor-reform}  Let $B$ be a scheme and let 
$f:D~\rightarrow~ B$ be a flat morphism with $D$ Gorenstein.  Let
$D_1\subset D$ be a codimension zero closed subscheme which is
Cohen-Macaulay.  Let $D_2\subset D$ be the residual
subscheme to $D_1\subset D$.  
\begin{enumerate}
\item For any closed subscheme $C\subset B$ which is a regular
embedding, $D_1\times_B C\subset D\times_B C$ and $D_2\times_B
C\subset D\times_B C$ are residual to each other.
\item If $B$ is regular, then $D_1$ and $D_2$ are flat over $B$.
\end{enumerate}
\end{cor}

\subsection{Reminder about Abel-Jacobi Maps}

We shall make occasional use of the \emph{Abel-Jacobi maps} associated 
to families of 1-cycles on $X$.  The reader is referred to
~\cite{CG},~\cite{VHS} for full definitions.  Here we recall only a
few facts 
about Abel-Jacobi maps.

\ps

Associated to a smooth, projective threefold $X$ there is a complex
torus  
\begin{equation}
J^2(X) = H^3_\ZZ(X)\backslash H^3(X,\CC) / \lt(H^{3,0}(X)\oplus
H^{2,1}(X)\rt).
\end{equation}
In case $X$ is a cubic hypersurface in $\PP^4$ (in fact for any
\emph{rationally connected} threefold) then $J^2(X)$ is a principally
polarized
abelian variety with theta divisor $\Theta$.  Given an algebraic
1-cycle $\gamma\in A_1(X)$ which 
is \emph{homologically equivalent to zero}~\cite[13]{VHS}, one can
associate a point $u_2(\alpha)$.  The construction is analogous to
the Abel-Jacobi map for a smooth, projective algebraic curve $C$ which 
associates to each 0-cycle $\gamma\in A_0(C)$ which is homologically
equivalent to zero a point $u_1(\alpha)\in J^1(C)$, the Jacobian
variety of $C$.  In particular $u_2:A_1(X)^{hom}\rightarrow J^2(X)$ is 
a group homomorphism.

\ps

Suppose that $B$ is a normal, connected variety of dimension $n$ and
$\Gamma\in A_{n+1}(B\times X)$ is an $(n+1)$-cycle such that for each
closed point $b\in B$ the corresponding cycle $\Gamma_b\in
A_1(X)$~\cite[\S 10.1]{F} is homologically equivalent to zero.  Then
in this case the set map $b\mapsto u_2(\Gamma_b)\in J^2(X)$ comes from 
a (unique) algebraic morphism $u=u_\Gamma:B\rightarrow J^2(X)$.  We
call this morphism the \emph{Abel-Jacobi map} determined by $\Gamma$.

\ps

More generally, suppose $B$ as above, $\Gamma\in A_{n+1}(B\times
X)$ is any $(n+1)$-cycle, and suppose $b_0\in B$ is some base-point.
Then we can form a new cycle $\Gamma'=\Gamma - \pi_2^*\Gamma_{b_0}$,
and for all $b\in B$ we have $\Gamma'_b = \Gamma_b - \Gamma_{b_0}$ is
homologically equivalent to zero.  Thus we have an algebraic morphism
$u=u_{\Gamma'}:B\rightarrow J^2(X)$.  Of course this morphism depends
on the choice of a base-point, but changing the base-point only
changes the morphism by a constant translation.  Thus we shall speak
of any of the morphisms $u_{\Gamma'}$ determined by $\Gamma$ and the
choice of a base-point as an \emph{Abel-Jacobi map} determined by
$\Gamma$.  

\ps

Suppose that $\Gamma_1,\Gamma_2\in A_{n+1}(B\times X)$ are two
$(n+1)$-cycles.  Then $u_{\Gamma_1+\Gamma_2}$ is the pointwise sum
$u_{\Gamma_1}+u_{\Gamma_2}$.  This trivial observation is frequently
useful.  Another useful observation is that any Abel-Jacobi morphism
$\alpha_\Gamma$ contracts all rational curves on $X$, since an Abelian 
variety contains no rational curves.

\section{Lines, Conics and Plane Cubics}~\label{sec-123}

We begin our analysis of the spaces $\Hdg{d}{g}(X)$ by recalling known 
results about the Fano scheme of lines on $X$, $F:=\Hdg{1}{0}(X)$.  

\ps

Two general lines $L_1,L_2\subset \PP^4$ determine a hyperplane by
$\text{span}(L_1,L_2)$.  We generalize this as follows:
Let $(F\times F-\Delta)\xrightarrow{\Phi} 
\PP^{4\vee}$ denote the following set map:
\begin{equation}\label{eqn-16}
 \Phi \left( [L_1,L_2] \right) =
  \left\{ \begin{array}{ll}
    \lt[ \text{span} \left( L_1,L_2 \right)\rt] &\text{if }L_1\cap L_2
  = \emptyset,  \\
  \lt[T_p X\rt]  & \text{if }p\in L_1\cap L_2
 \end{array} \right.
\end{equation} 
By~\cite[lemma 12.16]{CG}, $\Phi$ is algebraic.  
Let $X^\vee\subset \PP^{4\vee}$ denote the dual variety of $X$,
i.e., the variety 
parametrizing tangent hyperplanes to $X$.  Let $X^\vee_s\subset
X^\vee$ denote 
the subvariety parametrizing hyperplanes $H$ which are tangent to $X$ and
such that the singular locus of $H\cap X$ is not simply a single ordinary 
double point.  Let $U_s\subset U\subset F\times F$ denote the open
sets $\Phi^{-1}\lt(\PP^{4\vee} - X^\vee\rt)\subset \Phi^{-1}\lt(\PP^{4\vee} - 
X^\vee_s\rt)$.  Finally, let $I\subset F\times F$ denote the divisor
parametrizing 
incident lines, i.e., $I$ is the closure of the set $\{\lt([L_1],[L_2]\rt) 
: L_1\neq L_2, L_1\cap L_2\neq\emptyset\}$.
In~\cite{CG}, Clemens and Griffiths completely describe both
the total 
Abel-Jacobi map $F\times F\xrightarrow{\psi} J(X)$ and the
Abel-Jacobi map $F\xrightarrow{i} J(X)$.  Here is a summary of their results

\ps

\begin{thm} \label{thm-cg} \ \\
\begin{enumerate} 

\item The Fano variety $F$ is a smooth surface and the 
Abel-Jacobi
map $F\xrightarrow{u} J(X)$ is a closed immersion~\cite[theorem 7.8,
theorem 12.37]{CG}. 
\extraitemsep

\item The induced map $Alb(F)=J^2(F)~\rightarrow~ J(X)$ is an isomorphism of 
principally polarized Abelian varieties~\cite[theorem 11.19]{CG}.
\extraitemsep

\item The class of $u(F)$ in $J(X)$
is $\frac{[\Theta]^3}{3!}$~\cite[proposition 13.1]{CG}.
\extraitemsep

\item The difference of 
Abel-Jacobi maps 
\begin{equation}
\psi:F\times F \rightarrow J(X), \ \ \psi([L],[L']) = u([L]) - u([L'])
\end{equation}
maps $F\times F$ generically $6$-to-$1$ to the theta divisor
$\Theta~\subset ~J(X)$~\cite[section 13]{CG}.
\extraitemsep

\item Let $(\Theta-\{0\}) \xrightarrow{\mathcal{G}} 
\PP(H^{1,2}(X)^\vee)$ denote the Gauss map.  If we identify
$\PP(H^{1,2}(X))$ with $\PP^4$ via the
Griffiths residue calculus~\cite{GR}, then the composite map
\begin{equation}
(F\times_C F-\Delta)\xrightarrow{\psi}(\Theta -\{0\})\xrightarrow{\mathcal{G}}
\PP^{4\vee}
\end{equation}
is just the map $\Phi$ defined above~\cite[formula 13.6]{CG}.
\extraitemsep

\item The fibers of the Abel-Jacobi map form a \emph{ Schl\"afli
double-six}, i.e., the general fiber of $\psi:F\times F~\rightarrow~ J$
is of the form $\{(E_1,G_1),\dots,(E_6,G_6)\}$ where the lines
$E_i,G_j$ lie in a smooth hyperplane section of $X$, the $E_i$ are
pairwise skew, the $G_j$ are pairwise skew, and $E_i$ and $G_j$ are
skew iff $i=j$.  

There is a more precise result than above.  Let 
\begin{equation}R' \subset \lt(U\times_{\PP^{4\vee}} U\rt) \times F \times
\GR{}{3}{V}\times \GR{}{3}{V}\end{equation} 
be the closed subscheme parametrizing data
$\lt(\lt([L_1],[L_2]\rt),\lt([L_3],[L_4]\rt),[l], [H_1], [H_3]\rt)$
such that for each $i=1,\dots,4$, $l\cap L_i\neq\emptyset$ and such
that $H_1\cap X = l\cup L_1\cup L_4$, $H_2\cap X = l\cup L_2\cup L_3$.
Let $R\subset U\times_{\PP^{4\vee}} U$ be the image of $R'$ under the
projection map.  Let $\Delta\subset U\times U$ be the diagonal.  Then
the fiber product $U\times_\Theta U\subset U\times U$ is just the
union $R\cup\Delta$~\cite[p. 347-348]{CG}
\extraitemsep

\item The branch locus of $\Theta\xrightarrow{\mathcal{G}} \PP^{4\vee}$ equals
the branch locus of $F\times F \xrightarrow{\Phi} \PP^{4\vee}$ equals   
the dual variety of $X$, i.e., the variety 
parametrizing the tangent hyperplanes to $X$.  The 
ramification locus of $U \xrightarrow{\Phi} \PP^{4\vee}$ equals the
ramification 
locus of $U\xrightarrow{\psi} \Theta$ equals
the divisor $I$.  Each
such pair is a simple ramification point of both $\psi$ and
$\Phi$~\cite[lemma 13,8]{CG}.  
\end{enumerate}
\end{thm}

\subsection{Conics}\label{sec-cons}
Next we consider $\Hdg{2}{0}(X)$ which parametrizes plane conics on
$X$.  We are mostly interested just in the irreducibility of the spaces
$\Hdg{d}{g}(X)$, but in this case we can give a complete description
of $\Hdg{2}{0}(X)$.  We begin by proving that $\Hdg{2}{0}(X)$ is
smooth.

\ps

\begin{lem}~\label{lem-cons1}
$\Hdg{2}{0}(X)$ is smooth of dimension $4$.
\end{lem}

\begin{proof}

Any plane conic $C$ is a local complete
intersection.  So by lemma~\ref{lem-def1}, it suffices to prove that
$h^1(N_{C/X})=0$.  In fact we will prove that for each smooth conic
$C\subset X$, either $N_{C/X}~\cong~\OO_C(1)~\oplus~\OO_C(1)$ or else
$N_{C/S}~\cong~\OO_C~\oplus~\OO_C(2)$.    

\ps

We have the standard normal bundle
sequence:
\begin{equation}
\begin{CD}
 0 @>>> N_{C/X} @>>> N_{C/\PP^4} @>>> N_{X/\PP^4}|_C @>>> 0.
\end{CD} 
\end{equation}
Of course $N_{X/\PP^4}\cong \OO_{\PP^4}(3)|_C$ and it isn't hard to
see that 
\begin{equation}
N_{C/\PP^4}\cong \OO_{\PP^4}(2)|_C \oplus
\OO_{\PP^4}(1)^2|_C \cong \OO_C(4)\oplus \OO_C(2)\oplus \OO_C(2).
\end{equation}

\ps

By the Lefschetz hyperplane theorem we know that the 2-plane
$P=\text{span}(C)$ is not contained in $X$.  Therefore the induced map
$N_{C/P} ~\rightarrow~ N_{X/\PP^4}$ is injective with length $2$
cokernel.  It follows then that $N_{C/X}$, considered as a subsheaf of 
$N_{C/\PP^4}$ maps injectively to the quotient $N_{P/\PP^4}|_C \cong
\OO_C(2)\oplus \OO_C(2)$ and the cokernel is the length two cokernel
above.  So $N_{C/X}$ has degree $2$ and no summand of $N_{C/X}$ can
have degree higher than 2.  So either $N_{C/X}\cong \OO_C(1)\oplus
\OO_C(1)$ or else $N_{C/X}\cong \OO_C\oplus \OO_C(2)$.

\ps

Since $h^1(\OO_C(d)) = 0$ for all $d> -2$, we conclude that
$\Hdg{2}{0}(X)$ is smooth.
\end{proof}

Every plane conic $C\subset \PP^4$ is contained in a unique
$2$-plane $\text{span}(C)\subset \PP^4$.  Therefore over
$\Hdg{2}{0}(X)$ we have a flat family of $2$-planes, $\Pi\subset
\Hdg{2}{0}(X)\times \PP^4$ such that $\Pi_{[C]} = \text{span}(C)$.  Of
course the projection morphism $\Pi~\rightarrow~\Hdg{2}{0}(X)$ is
smooth.  By lemma~\ref{lem-cons1}, it follows that $\Pi$ is smooth.
Now consider the intersection $D\subset \Hdg{2}{0}(X)\times X$ of
$\Pi$ with $\Hdg{2}{0}(X)\times X$ in  
$\Hdg{2}{0}(X)\times \PP^4$.  First of all note that
$D~\rightarrow~\Hdg{2}{0}(X)$ 
has constant fiber dimension 1 over $\Hdg{2}{0}(X)$, since by the
Lefschetz hyperplane theorem~\cite[p. 156]{GH} $X$ contains no
$2$-planes.  Since 
$\Hdg{2}{0}(X)\times X$ is a Cartier divisor in $\Hdg{2}{0}(X)\times
\PP^4$, also $D\subset \Pi$ is a Cartier divisor.  In particular $D$
is a local complete intersection.  Therefore $D\rightarrow
\Hdg{2}{0}(X)$ is flat.  

\ps

Now let $\mathcal{C}\subset \Hdg{2}{0}(X)\times X$ denote the
universal smooth family of plane conics.  Then $\mathcal{C}$ is
smooth and $\mathcal{C}\subset D$ is a codimension zero closed
subscheme.  Let $D_2\subset D$ be the residual to $\mathcal{C}$ in
$D$.  Then by corollary~\ref{cor-reform}, we conclude that
$D_2\rightarrow \Hdg{2}{0}(X)$ is flat and the fiber of $D_2$ over a
closed 
point $[C]\in\Hdg{2}{0}(X)$ is simply the residual of $C$ in
$\text{span}(C)\cap X$.  But $\text{span}(C)\cap X$ is a plane cubic
curve, so the fiber of $D_2$ is just a line.  So we have an induced
morphism $g:\Hdg{2}{0}(X)~\rightarrow~ F$ which associates to
each $[C]$ the residual line in $\text{span}(C)\cap X$.  

\ps

Define $Q$ to be the rank 3 vector bundle on $F$ which is the quotient 
of $\OO_F^5$ by the universal sub-bundle.  
Let $\pi:\PP(Q)~\rightarrow~F$ be the projective bundle associated to
the rank 3 vector bundle.  The points of $\PP(Q)$ correspond to pairs
$([L],[P])$ where $L\subset X$ is a line and $P\subset \PP^4$ is a
$2$-plane such that $L\subset P$.  Therefore over $\PP(Q)$ we have a
flat family of $2$-planes $\Pi'\subset \PP(Q)\times \PP^4$.  
Let $D'\subset \Pi'$ denote the intersection of $\Pi'$ with $\PP(Q)
\times X$ and
let
$\mathcal{C}'\subset \PP(Q)\times X$ denote the pullback from $F$ of
the universal family of lines.  Then, $\mathcal{C}'\subset D'$ and the 
residual $D_2'$ is a flat family of conics.  Thus there is an induced
morphism $h:\PP(Q)\rightarrow \HI{2t+1}{X}$.  It is easy to see that
$h$ is a bijection of closed points over the open subset
$\Hdg{2}{0}(X)\subset \HI{2t+1}{X}$.  Since both $\PP(Q)$ and
$\Hdg{2}{0}(X)$ are smooth, it follows by Zariski's main
theorem~\cite[p. 288-289]{RB} that $\Hdg{2}{0}(X)$ is isomorphic to an 
open subset of 
$\PP(Q)$ and $g$ corresponds to the projection morphism $\PP
Q\rightarrow F$.  

\ps

But we can say more:  since the Abel-Jacobi morphism 
$u:F\rightarrow J(X)$ is an embedding, $F$ contains no rational
curves.  Thus, all the rational curves in $\PP(Q)$ lie in fibers.
Since $h$ is finite over $\Hdg{2}{0}(X)$,
no fiber of $\PP(Q)\rightarrow F$ is contracted by $h$, thus no
rational curve in $\PP(Q)$ is contracted by $h$ (since all rational
curves in $\PP(Q)$ are numerically equivalent, if one is contracted they 
all are).
But by~\cite[theorem VI.1.2]{K}, the exceptional locus of $h$ is
ruled.  Thus we conclude that $h$ is a finite morphism.  It follows by 
Zariski's main theorem that $h:\PP(Q)\rightarrow \HI{2t+1}{X}$ is the
normalization of $\HI{2t+1}{X}$.  We summarize the results as follows:

\ps

\begin{prop}~\label{prop-cons1}  The morphism
$\Hdg{2}{0}(X)\rightarrow F$ is isomorphic to an open subset of a
$\PP^2$-bundle $\PP(Q)\rightarrow F$.  In particular, $\Hdg{2}{0}(X)$
is smooth and connected of dimension $4$.  Moreover $\PP(Q)$ is the
normalization of $\HI{2t+1}{X}$.  
\end{prop}

\subsection{Plane Cubics}\label{sec-pcubs}

Every curve $C\subset \PP^4$ with Hilbert polynomial $3t$ is a plane
cubic, and the 2-plane $P=\text{span}(C)$ is unique; we have
that $C=X\cap P$.  Therefore the Hilbert scheme $\HI{3t}{X}$ is just
the Grassmannian $\mathbb{G}(2,4)$ of 2-planes in $\PP^4$ and
$\Hdg{3}{1}(X)$ is just an open subset of $\mathbb{G}(2,4)$.    

\section{Twisted Cubics}~\label{sec-30}

In this section we prove the irreducibility of $\Hdg{3}{0}(X)$. 
But first we prove an enumerative result about the number of
$2$-secant lines to a curve $C\subset X$.

\ps

Given a smooth curve $C\subset X$ we want to consider the set of
$2$-secant lines to $C$ which lie in $X$.

\begin{defn}
For a smooth curve $C\subset X$ we define $B_C\subset F$ to be the
scheme parametrizing lines in $X$ which intersect $C$ in a scheme of
degree $2$ or more.
\end{defn}

A dimension count leads one to expect that $B_C$ is a $0$-dimensional
scheme.  What is the degree of this scheme?

\begin{lem}\label{lem-enum}  Suppose that $C\subset X$ is a smooth
curve of genus $g$ 
and degree $d$.  
Define $b(C)=\frac{5d(d-3)}{2}+6-6g$.  If
$B_C$ is not positive dimensional and if $b(C)\geq 0$, then the degree of
$B_C$ is $b(C)$.  
\end{lem}

\begin{proof}  This is a standard Chern class
argument.  We work in the Chow ring of $C\times C$.  Let
$\omega\in A^*(C)$ denote the first Chern class of $\OO_{\PP^4}(1)|_C$ so that
$\omega$ is algebraically equivalent to $d$ times the class of a
point.  Let $\omega_1,\omega_2\in A^*(C\times C)$ denote the
pullbacks of $\omega$ by the two projection maps.  Let
$C\xrightarrow{\Delta} C\times C$ denote the diagonal morphism.
Also let $\Delta$, $\Delta_*\omega\in A^*(C\times C)$ denote the
class of the image of $\Delta$ and the class of the pushforward by
$\Delta$ of $\omega$ respectively.  

\ps

Let $V$ be the underlying vector space of $\PP^4$, and 
$A_C\subset \GR{\CC}{2}{V}$ be the scheme parametrizing chords
to $C$ in $\PP^4$.  We adopt the following convention: for $p\in C$ we denote by
$\text{span}(p,p)$ the tangent line to $C$ at $p$.  
Then
we have a morphism $C\times C\xrightarrow{f} A_C$ by
$(p,q)\mapsto[\text{span}(p,q)]$.  Let $S$ be the universal
rank two subbundle of $V\otimes \OO_{C\times C}$ whose fibre over a point
$(p,q)$ corresponds to the line $\text{span}(p,q)$.
The inclusion $S~\rightarrow~
V\otimes_\CC\OO_ {C\times C}$ induces a morphism of schemes
$P:=\PP(S)~\rightarrow~ (C\times C)\times \PP^4$.
We have two sections of $P$ determined by $(p,q)\mapsto
p\in\text{span}(p,q)$ and $(p,q)\mapsto q\in\text{span}(p,q)$.  Let
$\mathcal{I}_1$ and $\mathcal{I}_2$ denote the ideal sheaves of these sections in $P$.
Since both of these sections are divisors, the ideal sheaf of their
scheme theoretic union is just $\mathcal{I}_1\cdot \mathcal{I}_2\cong \mathcal{I}_1\otimes_{\OO_P}
\mathcal{I}_2$.  
Let $g:P\rightarrow \PP^4$ be the inclusion of $P$ into 
$(C\times C)\times \PP^4$ followed by projection onto $\PP^4$.
Let $D$ be the set of points of $P$ which are sent into $X$ under this map,
with ideal sheaf $\mathcal{I}_D=g^{*}\mathcal{I}_X$. 
 The two sections are two subvarieties of $D$, and  therefore we have that $\mathcal{I}_D\hookrightarrow
\OO_P$ factors through the subsheaf $\mathcal{I}_1\cdot \mathcal{I}_2\hookrightarrow
\OO_P$, i.e., we have $\mathcal{I}_D\hookrightarrow \mathcal{I}_1\cdot \mathcal{I}_2$.  The ideal
sheaf of the residual to these sections inside of $D$ is just what we
obtain when we twist this last map, namely
$\mathcal{I}_D\otimes_{\OO_P}(\mathcal{I}_1\cdot \mathcal{I}_2)^\vee\hookrightarrow \OO_P$.  We wish
to determine when this residual subscheme contains fibers of the projection map
$P\xrightarrow {\pi}C\times C$.  Let us assume that a general chord to 
$C$ does not lie in $X$.  
Then $\mathcal{I}_D$ is isomorphic to the locally
free sheaf $\OO_{S}(-3)$.  So we may twist our inclusion to get
$\OO_P~\rightarrow~ (\mathcal{I}_D)^\vee\otimes_{\OO_P}\mathcal{I}_1\otimes_{\OO_P}\mathcal{I}_2$.
The pushforward of this map yields a map 
\begin{equation}
\OO_{C\times
C}\xrightarrow{\phi}\pi_*\lt({(\mathcal{I}_D)^\vee\otimes_{\OO_P}\mathcal{I}_1 \otimes_{\OO_P}
\mathcal{I}_2\rule{0cm}{0.35cm}}\rt).
\end{equation}
It is clear that the fiber $\pi^{-1}(p,q)$ will be contained in
$D$ iff the image of the constant section $1$ under this map vanishes
at the stalk of $(p,q)$.  Therefore we conclude that the fiber product
$(C\times C)\times_{A_C} B_C$ is precisely the zero scheme of
$\phi$.  One sees that $(\mathcal{I}_D)^\vee\otimes_{\OO_P}\mathcal{I}_1\otimes_{\OO_P}
\mathcal{I}_2$ is a locally free sheaf of fiber degree $1$, in particular it is
relatively ample.  Therefore the pushforward
$E:=\pi_*((\mathcal{I}_D)^\vee\otimes_{\OO_P}\mathcal{I}_1\otimes_{\OO_P} \mathcal{I}_2)$ is a
locally free sheaf of rank $2$.  So, if the zero locus of $\phi$ is
zero dimensional, then we see that the class of this locus in
$A^*(C\times C)$ is just $c_2(E)$.  So we are reduced to a Chern
class calculation.

\ps

What is the Chern class of $S$?  The two sections in the last
paragraph yield a map of locally free sheaves $pr_1^*(\OO_{\PP^4}(-1)|_C)\oplus
pr_2^*(\OO_{\PP^4}(-1)|_C) ~\rightarrow~ S$.  This is an injective map and the
cokernel is supported on the diagonal.  Using the fact that the
cokernel of $S$ in $V\otimes_\CC \OO_{C\times C}$ is locally free
and a simple snake lemma argument, one deduces that the cokernel
is isomorphic to the coherent sheaf $\OO_{C\times
C}(\Delta)\otimes_{\OO_{C\times C}}\Delta_*(\OO_{\PP^4}(-1)|_C)$.  So we
deduce that the Chern class of $S$ is $1-\omega_1-\omega_2+\Delta
+\omega_1\cdot\omega_2-\Delta_*\omega$.  Let $\eta$ denote the first
Chern class of $\OO_S(1)$.  One has exact sequences
\begin{equation}
\begin{CD}
0@>>>\OO_S(1)\otimes_{\OO_P}(pr_i\circ \pi)^*(\OO_{\PP^4}(-1)|_C) @>>>
	\OO_S(1)\otimes_{\OO_P}\pi^*(S) @>>> \mathcal{I}_i^\vee @>>> 0
\end{CD}
\end{equation}
for $i=1,2$.  Thus one deduces that the Chern classes of $\mathcal{I}_1$ and $\mathcal{I}_2$
are $1-\eta+\omega_2-\Delta$ and $1-\eta+\omega_1-\Delta$ respectively.
Of course the Chern class of $\mathcal{I}_D$ is simply $1-3\eta$.  Since $\mathcal{I}_D^\vee
\otimes_{\OO_P} \mathcal{I}_1\otimes_{\OO_P} \mathcal{I}_2$ is relatively ample, its higher
direct images vanish.  Thus we may calculate the second Chern class of $E$
by a simple application of the Grothendieck-Riemann-Roch theorem~\cite{F}. 
It turns out
to be $5\omega_1\cdot\omega_2-15\Delta_*\omega+6\Delta\cdot\Delta$.  If we 
work modulo algebraic equivalence and omitting the phrase ``class of a point'',
we have $\omega_1\cdot\omega_2=d^2$, $\Delta_*\omega=d$ and $\Delta\cdot\Delta
=\chi(C)=2-2g$.  Using the fact that the map $f$ is generically $2$-to-$1$,
we deduce that the degree of $B_C$ is $\frac{5d(d-3)}{2}+6-6g$.
\end{proof}

\ps

\begin{lem}~\label{lem-30}  
$\Hdg{3}{0}(X)$ is smooth of dimension $6$.
\end{lem}

\begin{proof}
By lemma~\ref{lem-def1} we need to prove that $h^1(N_{C/X})=0$ for all 
$[C]\in\Hdg{3}{0}(X)$.  Consider the normal bundle sequence
\begin{equation}
\begin{CD}
0 @>>> N_{C/X} @>>> N_{C/\PP^4} @>>> N_{X/\PP^4}|_C @>>> 0.
\end{CD}
\end{equation}
Of course for any twisted cubic $C$, we have that $H=\text{span}(C)$
is a hyperplane, and $N_{C/H}~\cong\OO_C(5)^2$.  Thus we conclude
that $N_{C/\PP^4}~\cong~\OO_C(5)^2\oplus\OO_C(3)$,  and
$N_{X/\PP^4}|_C~\cong~ \OO_C(9)$.  So $N_{C/X}$ is a rank 2 vector
bundle of degree $4$.  By Grothendieck's lemma about vector bundles on 
$\PP^1$, we conclude $N_{C/X}~\cong~\OO_C(a)\oplus\OO_C(4-a)$ for some 
$a\geq 2$.  But since $N_{C/X}$ is a subbundle of
$\OO_C(5)^2\oplus\OO_C(3)$, we conclude that $a\leq 5$.  In all four
cases $a=2,3,4$ and $5$, we see that $4-a> -2$ so that
$h^1(N_{C/X})=0$.
\end{proof}

Define 
\begin{equation}
I=I_{3,0}\subset \Hdg{3}{0}(X)\times F
\end{equation}
to be the closed subset parametrizing pairs $([C],[L])$ where $L$ is a
$2$-secant line to 
$C$, and define 
\begin{equation}
f=f_{3,0}:I\rightarrow \Hdg{3}{0}(X)
\end{equation}
to be the projection.  
By lemma~\ref{lem-enum}, we know that $f_{3,0}$ is surjective.
Notice also that none of the lines $L$ is  
a $3$-secant line, because any $3$ points on a twisted cubic are
linearly independent.  

\ps

Now given $([C],[L])\in I$, the reducible curve $C\cup L$ lies on a
pencil of quadric surfaces in the $3$-plane $P=\text{span}(C)$, and
the general member of this pencil is smooth.  Let
$J~\subset~ I\times \Hi{t^2+2t+1}(\PP^4)$ denote the locally closed
subset parametrizing triples $([C],[L],[S])$ where $S$ is a smooth
quadric surface containing $C\cup L$.  Then $J~\rightarrow~ I$ is
birational to a $\PP^1$-bundle, in particular given an irreducible
component $J_i$ of $J$ with image $I_i\subset I$, we have
$\dim (I_i)~=~\dim (J_i)-1$.  
By the Lefschetz hyperplane
theorem, $X$ does not contain the surface $S$, thus
$S\cap X\subset S$ is a Cartier divisor of type $(3,3)$ on $S$.  The
residual to $C\cup L\subset S\cap X$ is a divisor of type $(1,1)$ on
$S$, i.e., a conic $D\subset S$.  

\ps

\begin{thm}\label{thm-irr30} The space $\Hdg{3}{0}(X)$ is a smooth, 
irreducible $6$-dimensional variety.
\end{thm}

\begin{proof}
By lemma~\ref{lem-30}, every irreducible component of $\Hdg{3}{0}(X)$
has dimension $6$.
We will prove that there is a unique irreducible component of 
$I$ of dimension $d\geq 6$. 
Since $I~\rightarrow~ \Hdg{3}{0}(X)$ is surjective, 
this implies that $\Hdg{3}{0}(X)$ is irreducible.
In order to show this, 
we will prove that $J$ has a unique
irreducible component of dimension $7$.

\ps

We stratify $J$ into locally closed subsets $J_1$, $J_2$, according to
the type of 
the residual curve $D$.  If $D$ is a smooth conic, we say that $D$ is
the first type.  If $D$ is a reducible conic, we say that $D$ is the
second type.  Notice that $D$ cannot be a \emph{double line} because
it is a divisor of type $(1,1)$ on a smooth quadric surface.

\ps

\textbf{Second Type:}
First consider $J_2$ parametrizing triples $([C],[L],[S])$ such that
$D$ is the second type.  Let $H\subset F\times F\times F$ denote the
locally closed subset parametrizing triples $([L],[D_1],[D_2])$ such
that $L$ and $D_1$ intersect transversally in one point, $D_1$ and
$D_2$ intersect transversally in one point, and $L$ is skew to $D_2$.
There is a morphism $J_2~\rightarrow~ H$ defined by decomposing
$D~=~D_1\cup D_2$ so that $L\cap D_1$ is non-empty.  Given a triple
$([L],[D_1],[D_2])$, every quadric surface $S$ containing $L\cup
D_1\cup D_2$ is contained in the $3$-plane $\text{span}(L,D_1,D_2)$.
Moreover, there is a $2$-dimensional linear system of quadrics $S$
containing $L\cup D_1\cup D_2$.  Thus the fiber dimension of
$J_2~\rightarrow~ H$ is at most $2$.  We can also see that the dimension
of $H$ is $4$:  there is a $2$-parameter family of choices for the
line $D_1$, and given $D_1$ there is a $1$-parameter family of lines
intersecting $D_1$.  Thus the dimension of $H$ is $2+1+1~=~4$.  So every 
irreducible component of $J_2$ has dimension at most $6$, which is
less than $7$.

\ps

Next we consider $J_1$ parametrizing triples $([C],[L],[S])$ such that 
the residual curve $D$ is a smooth conic.
Let $K\subset F\times \Hdg{2}{0}(X)$ denote the closed subset
parametrizing pairs $([L],[D])$ such that $L$ and $D$ intersect
transversally in one point $p$.  There is a morphism $J_1~\rightarrow~
K$ by sending $([C],[L],[S])$ to $([L],[D])$ with $D$ the residual
curve.
For a point $([L],[D])\in K$ and a point $([C],[L],[S])$ in the
fiber over $([L],[D])$, we have that $S$ is contained in the
$3$-plane $\text{span}(L,D)$.  There is a $2$-parameter
linear system of quadric surfaces $S\subset\text{span}(L,D)$ which contain
$L\cup D$.  The collection of quadric surfaces
$S\subset\text{span}(L,D)$ containing $L\cup D$ and such that also the 
residual curve $C$ of $L\cup D\subset S\cap X$ is a smooth twisted
cubic forms an open subset of the collection of all quadric surfaces
$S\subset \text{span}(L,D)$ containing $L\cup D$. 
So every (non-empty) fiber of $J_1~\rightarrow~ K$ is irreducible of
dimension $2$.  

\ps

Since $J_1~\rightarrow~ K$ has irreducible fibers of dimension $2$
(when they are non-empty), we see that for each irreducible component
$K_i$ of $K$, there is at most one
irreducible component of $J_1$ which fibers over $K_i$ with fiber
dimension $2$.  So we are reduced to proving that $K$ is irreducible
of dimension $5$.
In order to specify a pair $([L],[D])$ intersecting
at the point $p$, it is equivalent to specify $L$, a point $p\in L$,
and the line $N$ residual to $D$ since then $D$ is determined as the
conic residual to $N\subset X\cap\text{span}(N,p)$.  So $K$ is
isomorphic to an open subscheme of the product of the universal line
over $F$ (parametrizing pairs $(L,p)$) with another copy of $F$
(parametrizing $N$), and this is an irreducible
$5$-fold.  Thus there is at most one irreducible component of $J_1$ of 
dimension at least $7$, and such an irreducible component is exactly
$7$ dimensional.  All that remains is to show that at least one such 
component exists.

\ps

Since $\Hdg{3}{0}(X)$ is nonempty, and $J\rightarrow \Hdg{3}{0}(X)$ is
surjective with fibre dimension one, we conclude that $J_1$ has such a 
component, and therefore that 
$\Hdg{3}{0}(X)$ is an irreducible $6$-dimensional variety. 
\end{proof}

\subsection{The Abel-Jacobi map for $\Hdg{3}{0}(X)$}

In order to analyze $\Hdg{4}{0}(X)$ we will need to understand the
Abel-Jacobi map $u:\Hdg{3}{0}(X)\rightarrow J(X)$.

\ps

We have a morphism 
\begin{equation}
\Hdg{3}{0}(\PP^4)\xrightarrow{\sigma^{3,0}} \PP^{4\vee}
\end{equation}
defined by sending $[C]$ to $\text{span}(C)$.
This morphism makes $\Hdg{3}{0}(\PP^4)$ 
into a locally trivial bundle over $\PP^{4\vee}$ with fiber
$\Hdg{3}{0}(\PP^3)$ .  Recall from section~\ref{sec-123} that we defined
$X^\vee
\subset \PP^{4\vee}$ to be the dual variety of $X$ which parametrizes
tangent hyperplanes to $X$ and we defined $U$ to
be the complement of $X^\vee$ in $\PP^{4\vee}$.  Then we define 
$\Hdg{3}{0}_{U}(X)$ to be the
open subscheme of $\Hdg{3}{0}(X)$ which parametrizes twisted cubics,
$C$, in $X$ 
such that $\sigma^{3,0}([C])\in U$.  By the graph construction we may
consider $\Hdg{3}{0}_{U}(X)$ as a locally closed subvariety of
$U\times \HI{3t+1}{X}$.  
Let $\overline{\mathcal{H}}\subset U\times\HI{3t+1}{X}$ denote the
closure of $\Hdg{3}{0}_{U}(X)$ with the reduced induced scheme structure.
Denote by $\overline{\mathcal{H}} \xrightarrow{f} U$ the projection map. 

\ps

\begin{thm}\label{thm-st30}  Let
$\overline{\mathcal{H}}\xrightarrow{f''} U'\xrightarrow{f'} U$ 
be the Stein factorization of $\overline{\mathcal{H}}\xrightarrow{f} U$.
Then 
$\overline{\mathcal{H}}\xrightarrow{f''} U'$ is isomorphic to a $\PP^2$-bundle
$\PP_{U'}(E)\xrightarrow{} U'$ with $E$ a locally free sheaf of rank $3$,
and $U'\xrightarrow{f'} U$ is an unramified finite morphism of degree
$72$.  Moreover, the Abel-Jacobi map $\overline{\mathcal{H}}\xrightarrow{i}
J(X)$ 
factors as $\overline{\mathcal{H}}\xrightarrow{f''} U'\xrightarrow{i'}
J(X)$ where 
$U'\xrightarrow{i'} J(X)$ is a birational morphism of $U'$ to a translate
of $\Theta$.
\end{thm} 

\begin{proof}

We need to use the following lemma:

\begin{lem}\label{lem-cubicsurf}
Let $S$ be a smooth cubic surface in $\PP^3$. Then there are exactly $72$
line bundles $L$ on $S$ such that $L^2=1$, and $L.K_S=-3$ 
$($where $K_S$ is the canonical class$)$.  Furthermore, each of them
satisfies $H^1(S,L)=H^2(S,L)=0$, and the general member of $H^0(S,L)$ is
a smooth curve.
\end{lem}

We will explicitly describe such bundles $L$ below, and this lemma
will be a straightforward consequence.  Note that if $C\subset S$ is a curve
with Hilbert polynomial $3t+1$, then $C.K_S=-3$, since $K_S$ is 
minus the hyperplane class, and since the curve has arithmetic genus zero,
adjunction shows that $C^2=1$. This shows that all the curves 
in $\mathcal{H}_U^{3,0}(X)$ give line bundles $L$ satisfying the conditions 
above.  Conversely, given any effective divisor $C\in |L|$, with $L$ a line
bundle as above, we see that $C$ has degree three, and arithmetic genus zero,
and hence Hilbert polynomial $3t+1$.

\ps

Now, let $\mathcal{X}\xrightarrow{\pi} U$ be the universal family of smooth
hyperplane sections of $X$.  For any $[H]\in U$, we use $S_H:=H\cap X$ to denote
the smooth cubic surface which is the fibre of $\pi$. Let 
$\pic^{3,0}(\mathcal{X}/U)$ be the subscheme of the relative Picard scheme 
parameterizing line bundles $L_H$ on the fibers $S_H$ of $\pi$ such that
$L_H^2=1$ and $L_H\cdot K_{S_H}=-3$.  For any such $L_H$, we have
$\chi_{S_H}(L_H)=1/2(L_H^2-L\cdot K_{S_H})+\chi(\OO_{S_H})=2+1=3$.
By the above lemma, the line bundle $L_H$ has no higher cohomology on $S_H$,
and so there is a rank three vector bundle $E$ on $\pic^{3,0}(\mathcal{X}/U)$
whose fibre at a point $(H,L_H)$ of $\pic^{3,0}(\mathcal{X}/U)$ consists of
the global sections $H^0(S_H,L_H)$. Let $P=\PP(E)$, be the projectivization
of this bundle, with projection map 
$g:P\longrightarrow \pic^{3,0}(\overline{X}/U)$.  
A point of this projectivization consists of the data $(H,L_H,C_H)$ where
$[H]\in U$, $L_H$ is a line bundle on $S_H$ satisfying the numerical conditions,
and $C_H$ is an effective divisor on $S_H$ with $\OO_{S_H}(C_H)=L_H$.  

\ps
By the remarks after lemma \ref{lem-cubicsurf}, we see that $C_H$ has 
Hilbert polynomial
$3t+1$, and so we have a natural map $P\longrightarrow U\times \Hi{3t+1}(X)$
sending $([H],L_H,C_H)$ to $([H],C_H)$.  The map is clearly an injection,
since we can recover the line bundle $L_H$ from $C_H$.  

\ps
The short exact sequence
$$0\longrightarrow \OO_{S_H} \longrightarrow \OO_{S_H}(C_H) 
\longrightarrow \OO_{S_H}(C_H)|_{C_H} \longrightarrow 0$$
gives the long exact sequence in cohomology
$$ 0  
     \longrightarrow 
   H^0\left({S_H,\OO_{S_H}\rule{0cm}{0.35cm}}\right)
     \xrightarrow{\cdot C_H}
   H^0\left({S_H,\OO_{S_H}(C_H)\rule{0cm}{0.35cm}}\right)
     \longrightarrow
   H^0\left({C_H,N_{C_H/S_H}\rule{0cm}{0.35cm}}\right)
     \longrightarrow
   H^1\left({S_H,\OO_{S_H}\rule{0cm}{0.35cm}}\right) = 0.
$$
This sequence has the following interpretation.  $H^0(S_H,\OO_{S_H}(C_H))$
divided by the section $C_H$ is the vertical tangent space (at $C_H$) for
the map $P\longrightarrow \pic^{3,0}(\mathcal{X}/U)$.  $H^0(C_H,N_{C_H/S_H})$
is the tangent space (at $C_H$) of $\Hi{3t+1}(S_H)$.  The sequence above
shows that the map between these tangent spaces is an isomorphism, and 
hence that $P\longrightarrow U\times \Hi{3t+1}(X)$ is a closed embedding.

\ps

The subset $\mathcal{H}_U^{3,0}(X)$ of $U\times\Hi{3t+1}(X)$ is contained 
in the image
of $P$, and, by lemma \ref{lem-cubicsurf}, $\mathcal{H}_U^{3,0}(X)$ is dense
in each fibre of $P\longrightarrow U$, so we conclude that $P=
\overline{\mathcal{H}}$.
The map $\rho:\pic^{3,0}(\mathcal{X}/U)\longrightarrow U$, since it is a 
finite type subscheme of the the relative Picard scheme, and by 
lemma \ref{lem-cubicsurf} each fibre
consists of $72$ points, i.e., the map is finite.  We also know that this
is unramified since the Picard group of a cubic surface is reduced.  Finally,
the map $P\longrightarrow \pic^{3,0}(\mathcal{X}/U)$ is a $\PP^2$ bundle
by construction.

\ps
Therefore, $P\xrightarrow{g}\pic^{3,0}(\mathcal{X}/U)\xrightarrow{\rho}
U$ is the Stein factorization
$\overline{\mathcal{H}}\xrightarrow{f''} U'\xrightarrow{f'} U$ of
$\overline{\mathcal{H}}\xrightarrow{f} U$, and this factorization has the 
properties claimed in the theorem.

\ps

It only remains to determine the Abel-Jacobi map
$\overline{\mathcal{H}}\xrightarrow{i} 
J(X)$.  Since $J(X)$ contains no rational curves, $i$ is a
constant map on each fiber 
of $g$.  Since $\text{Pic}^{3,0}(\mathcal{X}/U)$ is smooth, it follows that $i$
factors through a morphism $i':\text{Pic}^{3,0}(\mathcal{X}/U)~\rightarrow~ J(X)$.  
To determine $i'$, we introduce the locus of ``Z's of lines'' i.e.
the subscheme of $\overline{\mathcal{H}}$ parametrizing cubic curves whose
irreducible components are lines.  To be precise, let $\Sigma\subset
\overline{\mathcal{H}}\times X$ be our
flat family of cubic curves.  We let $\Sigma^{s}\subset\Sigma$ denote the
singular subscheme.  We can form the flattening stratification for $\Sigma^{s}
~\rightarrow~ \overline{\mathcal{H}}$, and we define
$Z\subset\overline{\mathcal{H}}$ to be
the stratum corresponding to the constant Hilbert polynomial $2$,
i.e., the locus parametrizing curves with two nodes.  What are
the fibers $Z\cap g^{-1}(q)$ for $q\in \text{Pic}^{3,0}(\mathcal{X}/U)$?  
Define $H~=~\rho(q)$.  In the analysis below, we will see that we can
find a set of $6$ mutually skew lines in $S_H$ such that $g^{-1}(q)$
corresponds to the complete linear series of lines in the blown-down $\PP^2$.
It is clear that a line $\ell$ in this linear series will 
correspond to a singular
cubic curve iff $\ell$ intersects one of the $6$ special points.  Similarly, $\ell$ 
will correspond to a cubic curve whose singular locus has degree $2$ iff
$\ell$ is one of the $15$ lines joining a pair of the $6$ special points.
Thus each
fiber $Z\cap g^{-1}(q)$ consists of $15$ points, and  
$\Sigma_Z^{s}~\rightarrow~ 
Z$ is an unramified, finite morphism of degree $2$. Thus
$\Sigma_Z^s~\rightarrow~ \text{Pic}^{3,0}(\mathcal{X}/U)$ is an unramified, finite
morphism of degree $30$.  

\ps

Denote by $i_1:\Sigma^s_Z\rightarrow J(X)$ the composition of
$\Sigma^s_X\rightarrow \text{Pic}^{3,0}(\mathcal{X}/U)$ with the Abel-Jacobi map
$\text{Pic}^{3,0}(\mathcal{X}/U)\rightarrow J(X)$.
Recall $\Sigma_Z^{s}$ parametrizes pairs 
$\lt([C],[x]\rt)$, where $C$ is a completely reducible cubic and $x$ is 
a node of $C$.  We define a map $h:\Sigma^{s}~\rightarrow~ F\times F$ as
follows.  The union of those components of $C$ which intersect $x$ is a
completely reducible conic, $C'$.  The residual to $C'$ inside of $C$ is
a line $\ell_1$.  Now $C'$ spans a $\PP^2$ in $\PP^4$ and the residual to
$C'$ in $\text{span }(C')\cap X$ is another line $\ell_2$.  We define $h$
to be the map $\lt([C],[x]\rt)\mapsto\lt([\ell_1],[\ell_2]\rt)$.  The point is,
since $C'$ and $\ell_2$ are residual in a complete intersection which
varies in a rational family, it follows
by the residuation trick that $i_1$ is equal to $\psi\circ h$
(up to a fixed translation).  

\ps

What are the fibers of $h$?  Suppose we are
given two skew lines $\ell_1$ and $\ell_2$ whose span intersects $X$ in a smooth
cubic surface, $X'$.  How many reducible conics $C'$ are there which are
residual to $\ell_2$ and which intersect $\ell_1$?  One of the lines in $C'$,
call it $\ell_3$, intersects both $\ell_1$ and $\ell_2$.  The other line of
$C'$ is uniquely determined by the condition that it be residual to
$\ell_2\cup \ell_3$ in the $\PP^2$ they span.  Thus the points in a fiber 
of $h$ are 
enumerated by the lines $\ell_3$ joining $\ell_1$ and $\ell_2$.  There are $5$ such
lines.  Therefore $h$ is dominant and generically finite of degree $5$.  
We know that 
$\psi$ maps dominantly and generically finitely to $\Theta$ of degree $6$,
thus $\Sigma_Z^{s}$ maps to $\Theta$ dominantly and generically finitely of
degree $5\times 6~=~30$.  We have already seen that
$\Sigma_Z^s~\rightarrow~ \text{Pic}^{3,0}(\mathcal{X}/U)$ is unramified of degree $30$.
Therefore $\text{Pic}^{3,0}(\mathcal{X}/U)~\rightarrow~ J(X)$ maps generically 1-to-1
and dominates a translate of $\Theta$.
\end{proof} 

\ps

\begin{cor}\label{cor-aj30}  The Abel-Jacobi map
$i_{3,0}:\Hdg{3}{0}(X)~\rightarrow~ J(X)$ dominates a translate of
$\Theta$ and is birational to a $\PP^2$-bundle over its image.
\end{cor}

\begin{proof}  $\mathcal{H}_U^{3,0}(X)\subset \mathcal{H}^{3,0}(X)$ is dense 
in $\overline{\mathcal{H}}$.
\end{proof}

We now need to examine the line bundles $L$ on a cubic surface $S$
satisfying $L^2=1$ and $L.K_S=-3$.  We need to establish the 
facts claimed in lemma \ref{lem-cubicsurf}, and also show that for any such
$L$, we can always blow down six lines so that $L$ is pullback of
$\OO_{\PP^2}(1)$ from the resulting $\PP^2$.
We will follow~\cite{H}, chapter V, notation
4.7.3 for our notation of the Neron-Severi group of $S$.  Recall
that $e_1,\dots,e_6$ are the linear equivalence classes of $6$
mutually skew lines on $S$, so that the contraction of
$e_1,\dots,e_6$ is a $\PP^2$,  and $l$ is the linear equivalence class
of the total transform of a line in $\PP^2$.  If we write
$L~=~al-\sum b_ie_i$, then we have $3a-\sum b_i~=~3$ and $a^2-\sum
b_i^2 ~=~ 1$.  Since $(\sum b_i)^2\leq 6\sum b_i^2$, we deduce that
$(3a-3)^2\leq 6a^2-6$.  This implies that either $a~=~1,2,3,4$, or
$5$.  One quickly works out all the possibilities for
$(a,b_1,\dots,b_6)$.  There is an obvious action of the group $S_6$ on
the set of solutions via permuting $b_1,\dots,b_6$.  Representatives
of the orbits of the set of solutions are as follows:
\begin{equation} (1,0,0,0,0,0,0),\ \ (2,1,1,1,0,0,0),\ \ (3,2,1,1,1,1,0),\ \ 
(4,2,2,2,1,1,1),\ \ (5,2,2,2,2,2,2). \rule[-0.3cm]{0cm}{0.3cm} \end{equation}
Counting the size of each orbit shows that there are a total of $72$ distinct
solutions.  With a slight amount of work, one shows that the separate orbits 
all lie in the same orbit under the action of the full Weyl group of $E_6$.  
Thus, for some choice of $6$ mutually skew six lines, we have that $L$ is
just $l$.  The general member of this linear series is obviously smooth,
and the long exact sequence in cohomology coming from

$$0\longrightarrow \OO_{S}\longrightarrow \OO_{S}(l)\longrightarrow
\OO_{\PP^1}(1)\longrightarrow 0$$
shows that $H^1(S,L)=H^2(S,L)=0$, which was the last thing to be checked.

\section{Quartic Elliptic Curves}~\label{sec-41}

Recall that the normalization of $\HI{2t+1}{X}$ is isomorphic to the
$\PP^2$-bundle $\PP(Q)
~\rightarrow~ F$ which parametrizes pairs $(L,P)$ which $L\subset X$
a line and $P\subset\PP^4$ a $2$-plane containing $L$.    
Let $A\xrightarrow{g} \PP(Q)$ denote the $\PP^1$-bundle
which parametrizes triples $(L,P,H)$ with $H$ a hyperplane containing
$P$.  Let $I_{4,1}\xrightarrow{h} A$ denote the $\PP^4$-bundle parametrizing
4-tuples $(L,P,H,Q)$ where $Q\subset H$ is a quadric surface
containing the conic $C\subset X\cap P$. 
Notice that $I_{4,1}$ is smooth
and connected of dimension $4+1+4=9$.

\ps

Let $D\subset I_{4,1}\times X$
denote the intersection of the universal quadric surface over $I_{4,1}$ with 
$I_{4,1}\times X\subset I_{4,1}\times \PP^4$.  Then $D$ is a local complete
intersection scheme.  By the Lefschetz hyperplane theorem, $X$
contains no quadric surfaces; therefore $D~\rightarrow~ I_{4,1}$ has constant
fiber dimension $1$ and so is flat.  Let $D_1\subset I_{4,1}\times X$ denote
the pullback from $\PP(Q)\times X ~=~\HI{2t+1}{X}\times X$ of the
universal family of conics.  Since $I_{4,1}\times X~\rightarrow~ \PP(Q)\times
X$ is smooth and the universal family of conics is a local complete
intersection which is flat over $\PP(Q)$, we conclude that also $D_1$ is a
local complete intersection which is flat over $I_{4,1}$.  Clearly
$D_1\subset D$.  Thus by corollary~\ref{cor-reform}, we see that the
residual $D_2$ of $D_1\subset D$ is Cohen-Macaulay and flat over $I_{4,1}$.

\ps

By the base-change property in corollary~\ref{cor-reform}, we see that
the fiber of $D_2~\rightarrow~ I_{4,1}$ over a point $(L,P,H,Q)$ is simply the
residual of $C\subset Q\cap X$.  If we choose $Q$ to be a smooth
quadric, i.e., $Q\cong \PP^1\times \PP^1$, then $C\subset Q$ is a
divisor of type $(1,1)$ and $X\cap Q\subset Q$ is a divisor of type
$(3,3)$.  Thus the residual curve $E$ is a divisor of type $(2,2)$,
i.e., a quartic curve of arithmetic genus 1.  Thus $D_2\subset I_{4,1}\times
X$ is a family of connected, closed subschemes of $X$ with Hilbert
polynomial $4t$.  So we have an induced map $f:I_{4,1}~\rightarrow~
\HI{4t}{X}$.  

\ps

\begin{prop}\label{lem-qell}  The image of the morphism above
$f:I_{4,1}~\rightarrow~\HI{4t}{X}$ is the closure $\overline{\Hdg{4}{1}}(X)$ 
of $\Hdg{4}{1}(X)$.  Moreover the open set $f^{-1}\Hdg{4}{1}(X)\subset 
I_{4,1}$ is a $\PP^1$-bundle over $\Hdg{4}{1}(X)$.  Thus $\Hdg{4}{1}(X)$ is
smooth and connected of dimension 8.
\end{prop}

\begin{proof}
If $E\subset X$ is a smooth, connected curve with Hilbert polynomial
$4t$, then $E$ is a quartic elliptic curve in some hyperplane $H$.
Any such curve lies on a pencil of quadric surfaces $Q$,  and the
residual of $E\subset Q\cap X$ is a conic.  Thus we see that $f(I_{4,1})$
contains the open subscheme $\Hdg{4}{1}(X)\subset \HI{4t}{X}$.  Since
$f(I_{4,1})$ is closed and irreducible, we conclude that
$f(I_{4,1})~=~\overline{\Hdg{4}{1}}(X)$.  Since the fibre of $f$
over any smooth elliptic quintic $E$ is determined by the $\PP^1$ of 
quadrics $Q$ in $H=\text{span}(E)$, we see that
$f^{-1}(\Hdg{4}{1}(X))\subset I_{4,1}$ is an open subset which is a
$\PP^1$-bundle over $\Hdg{4}{1}$.  In particular, since $I_{4,1}$ is smooth
and connected,
we conclude that $\Hdg{4}{1}$ is also smooth and connected of
dimension 8.
\end{proof}

The surface $F$ of lines contains no rational curves, so in the $\PP^1$ fibre of
$f$ over $[E]\in\Hdg{4}{1}$, the line $L$ must be constant.  Since the 
hyperplane $H$ is also determined by $[E]$, we have a well-defined morphism
$m:\Hdg{4}{1}(X)~\rightarrow~ \PP(Q^\vee)$, where $\PP(Q^\vee)$ is the 
$\PP^2$-bundle over $F$ parametrizing pairs
$([L],[H]),\, L\subset H$.  For a general $H$, the intersection $Y~=~H\cap X$ is
a smooth cubic surface,  and the fiber $m^{-1}([L],[H])$ is an open
subset of the complete linear series $|\OO_Y{L+h}|$, where $h$ is the
hyperplane class on $Y$.  Thus $m:\Hdg{4}{1}(X)~\rightarrow~ \PP(Q^\vee)$
is a morphism of smooth connected varieties which is birational to a
$\PP^4$-bundle.  Composing $m$ with the projection $\PP(Q^\vee)$ yields
a morphism $n:\Hdg{4}{1}(X)~\rightarrow~ F$ which is birational to a
$\PP^4$-bundle over a $\PP^2$-bundle.  

\ps

\begin{cor}~\label{cor-41}
The morphism $n:\Hdg{4}{1}(X)~\rightarrow~ F$ from above is birational
over $F$ to a $\PP^4$-bundle over a $\PP^2$-bundle over $F$.
\end{cor}

\section{Cubic Scrolls and Applications}~\label{sec-cuscr}

\subsection{Preliminaries on cubic scrolls}~\label{subsec-prelim}

In the next few sections we will use residuation in a cubic 
scroll. We start by collecting 
some basic facts about these surfaces.

\ps

There are several equivalent descriptions of cubic scrolls.  
\begin{enumerate}
\item A cubic scroll $\Sigma \subset \PP^4$ is a connected, smooth
surface with Hilbert polynomial $P(t)  =  \frac{3}{2}t^2 + \frac{5}{2}t
+ 1$. 
\extraitemsep

\item A cubic scroll $\Sigma\subset
\PP^4$ is the determinantal variety defined by the $2\times 2$ minors of 
a matrix of linear forms:
\vskip -0.8cm
\begin{equation}\lt[ \begin{array} {ccc}
 L_1 & L_2 & L_3 \\
 M_1 & M_2 & M_3 
\end{array} \rt]\end{equation}
such that for each row or column, the linear forms in
that row or column are linearly independent 
\extraitemsep

\item A cubic scroll $\Sigma\subset \PP^4$ is the \emph{join} of an
isomorphism $\phi:L~\rightarrow~ C$.  Here $L\subset \PP^4$ is a line and
$C\subset \PP^4$ a conic such that $L\cap
\text{span}(C) ~=~\emptyset$.  The join of $\phi$ is defined as the
union over all $p\in L$ of the line $\text{span}(p,\phi(p))$.  
\extraitemsep

\item A cubic scroll $\Sigma\subset \PP^4$ is the image of a morphism $f:\PP(E)
\rightarrow~ \PP^4$ where 
$E~=~\OO_{\PP^1}(-1)\oplus \OO_{\PP^1}(-2)$ on $\PP^1$, the morphism $f:\PP(E)
~\rightarrow~ \PP^4$ is such that $f^*\OO_{\PP^4}(1) ~=~ \OO_{E}(1)$,
and the pullback map $H^0(\PP^4,\OO_{\PP^4}(1))\rightarrow H^0(\PP(E),
\OO_{E}(1))$ is an isomorphism.
\extraitemsep

\item A cubic scroll $\Sigma\subset \PP^4$ is as
a minimal variety, i.e., $\Sigma\subset \PP^4$ is any smooth connected 
surface with $\text{span}(\Sigma)~=~\PP^4$ which has the minimal
possible degree for such a surface, namely $\text{deg}(\Sigma)~=~3$. 
\extraitemsep

\item A cubic scroll $\Sigma\subset \PP^4$ is a smooth surface residual to a
2-plane $\Pi$ in the base locus of a pencil of quadric hypersurfaces
which contain $\Pi$.  
\end{enumerate}

\ps

From the fourth description $\Sigma ~=~ \PP(E)$ we see that
$\text{Pic}(\Sigma) ~=~ \text{Pic}(\PP(E)) \cong \ZZ^2$.  Let $\pi:\PP(E)
~\rightarrow~ \PP^1$ denote the projection morphism and let
$\sigma:\PP^1~\rightarrow~ \PP(E)$ denote the unique section whose image
$D~=~\sigma(\PP^1)$ has self-intersection $D.D~=~-1$.  Then $f(D)$ is a
line on $\Sigma$ called the \emph{directrix}.  For each $t\in
\PP^1$, $f(\pi^{-1}(t))$ is a line called a \emph{line of ruling} of
$\Sigma$.  Denote by $F$ the divisor class of any $\pi^{-1}(t)$.  Then
$\text{Pic}(\Sigma) ~=~ \ZZ\{D,F\}$ and the intersection pairing on
$\Sigma$ is determined by $D.D ~=~ -1, D.F ~=~1, F.F ~=~0$.  The hyperplane
class is $H~=~D+2F$ and the canonical class is $K ~=~ -2D-3F$.  

\ps

Using the fourth description of a cubic scroll, we see that any two
cubic scrolls differ only by the choice of the isomorphism
$H^0(\PP^4,\OO_{\PP^4}(1)) \rightarrow H^0(\PP(E),\OO_{E}(1))$.
Therefore any two cubic scrolls are
conjugate under the action of $\text{PGL}(5)$.  So the open set
$U\subset \text{Hilb}_{P(t)}(\PP^4)$ parametrizing cubic scrolls is a
homogeneous space for 
$\text{PGL}(5)$, in particular it is smooth and connected.

\ps

One possible specialization of a cubic scroll is a reducible surface
$\Sigma ~=~ \Sigma_1\cup \Sigma_2$ where $\Sigma_1$ is a 2-plane,
$\Sigma_2$ is a smooth quadric surface, and $\Sigma_1\cap \Sigma_2$ is
a line $L$.  Let $T\subset \text{Hilb}_{P(t)}(\PP^4)$ denote the locus
parametrizing surfaces $\Sigma$ of this form.  

\ps

\begin{lem}~\label{lem-specsm}  The Hilbert scheme
$\text{Hilb}_{P(t)}(\PP^4)$ is 
smooth along $T$.  
\end{lem}

\begin{proof}
Let $[\Sigma]$ be a point of $T$.  Since $\Sigma$ is a local complete
intersection, the normal sheaf $N_{\Sigma/\PP^4}$ is locally free.  The
Zariski tangent space of $\text{Hilb}_{P(t)}(\PP^4)$ at $[\Sigma]$ is
identified with $H^0(\Sigma,N_{\Sigma/\PP^4})$.  By \cite[theorem I.2.15.2]{K},
%[Koll\'ar96] theorem I.2.15.2, 
we see that every irreducible component of
$\text{Hilb}_{P(t)}(\PP^4)$ through $[\Sigma]$ has dimension at least 
\begin{equation}\dim  H^0(\Sigma,N_{\Sigma/\PP^4}) - \dim 
H^1(\Sigma,N_{\Sigma/\PP^4}).\end{equation} 
So once we show that $H^1(\Sigma,N_{\Sigma/\PP^4})~=~0$, it will follow
that $\text{Hilb}_{P(t)}(\PP^4)$ is smooth at $[\Sigma]$.

\ps

In order to analyze the normal bundle $N_{\Sigma/\PP^4}$ we recall the
following result: Suppose that $X$ is a smooth ambient
variety and suppose that $Y\subset X$ is a simple normal crossings
variety with no triple points.  Let $Y_i\subset Y$ be an irreducible
component and let $Z_1,\dots,Z_r$ be the connected components of
$\text{Sing}(Y)\cap Y_i$.  For each $i~=~1,\dots,r$, there is an \'etale
cover $f:W~\rightarrow~ V$ of a Zariski neighborhood of $Z_i\subset Y$
such that the preimage $Z~=~f^{-1}(Z_i)$ of $Z_i$ is connected and such that $W$
is reducible along $Z$.  Denote the two branches of $W$ along
$Z$ by $W'$ and $W''$.  Then the line bundle $N_{Z/W'}\otimes
N_{Z/W''}$ descends to a line bundle $N_i$ on $Z_i$.  We have a short
exact sequence of coherent sheaves:
\begin{equation}\begin{CD}
0 @>>> N_{Y_i/X} @>>> N_{Y/X}|_{Y_i} @>>> \oplus_{i~=~1}^r N_i @>>> 0.
\end{CD} \end{equation}
In our particular case, we have the two exact sequences:
\begin{equation}\begin{CD}
0 @>>> N_{\Sigma_1/\PP^4} @>>> N_{\Sigma/\PP^4}|_{\Sigma_1} @>>>
N_{L/\Sigma_1}\otimes N_{L/\Sigma_2} @>>> 0 \\
0 @>>> N_{\Sigma_2/\PP^4} @>>> N_{\Sigma/\PP^4}|_{\Sigma_2} @>>> 
N_{L/\Sigma_1}\otimes N_{L/\Sigma_2} @>>> 0.
\end{CD} \end{equation}

\ps

If we identify $\Sigma_1$ with $\PP^2$, then we have
$N_{\Sigma_1/\PP^4}\cong \OO_{\PP^2}(1)\oplus \OO_{\PP_2}(1)$.  If we
identify $\Sigma_1$ with $\PP^1\times\PP^1$, then we have
$N_{\Sigma_2/\PP^4} ~=~ \OO_{\PP^1\times \PP^1}(1,1) \oplus
\OO_{\PP^1\times\PP^1}(2,2)$.  Identifying $L$ with $\PP^1$, we have
$N_{L/\Sigma_1} \cong \OO_{\PP^1}(1)$ and $N_{L/\Sigma_2} \cong
\OO_{\PP^1}$.  We are more interested in
$N_{\Sigma/\PP^4}|_{\Sigma_2}(-L)$ than we are in
$N_{\Sigma/\PP^4}|_{\Sigma_2}$.  To relate the two we use the
identification $\OO_{\Sigma_2}(-L) \cong \OO_{\PP^1\times
\PP^1}(-1,0)$.  With all of these identifications, we get two exact
sequences: 
\begin{equation}\begin{CD}
0 @>>> \OO_{\PP^2}(1)\oplus \OO_{\PP^2}(1) @>>>
N_{\Sigma/\PP^4}|_{\Sigma_1} @>>> \OO_{\PP^1}(1) @>>> 0 \\
0 @>>> \OO_{\PP^1\times \PP^1}(0,1)\oplus \OO_{\PP^1\times \PP^1}(1,2)
@>>> N_{\Sigma/\PP^4}|_{\Sigma_2}(-L) @>>> \OO_{\PP^1}(1) @>>> 0.
\end{CD} \end{equation}
Applying the long exact sequence in cohomology to these two short
exact sequences, we conclude the vanishing result 
\begin{equation}H^1(\Sigma_1,N_{\Sigma/\PP^4}|_{\Sigma_1}) ~=~
H^2(\Sigma_1,N_{\Sigma/\PP^4}|_{\Sigma_1}) ~=~ \end{equation}
\begin{equation}H^1(\Sigma_2,N_{\Sigma/\PP^4}|_{\Sigma_2}(-L)) ~=~
H^2(\Sigma_2,N_{\Sigma/\PP^4}|_{\Sigma_2}(-L)) ~=~ 0.\end{equation}

\ps

We also have a short exact sequence:
\begin{equation}\begin{CD}
0 @>>> N_{\Sigma/\PP^4}|_{\Sigma_2}(-L) @>>> N_{\Sigma/\PP^4} @>>>
N_{\Sigma/\PP^4}|_{\Sigma_1} @>>> 0.
\end{CD}\end{equation}
Applying the long exact sequence in cohomology to this short exact
sequence and combining with the vanishing result of the last
paragraph, we conclude that $H^1(\Sigma,N_{\Sigma/\PP^4})  = 
H^2(\Sigma,N_{\Sigma/\PP^4}) = 0.$  Therefore
$\text{Hilb}_{P(t)}(\PP^4)$ is smooth along $T$.
\end{proof}

\ps

\begin{lem}~\label{lem-smscr}  
The union $V~=~T\cup U\subset \text{Hilb}_{P(t)}(\PP^4)$ is a
smooth, connected open subset.  
\end{lem}

\begin{proof}
Given $[\Sigma]\in T$, we will show that every deformation of $\Sigma$
can be realized as a subvariety of a rank 4 quadric hypersurface
$Q\subset\PP^4$.  Then we will examine the deformations of $\Sigma$ as
a subvariety of $Q$ to prove the lemma.

\ps

Let $I\subset \text{Hilb}_{P(t)}(\PP^4)\times \PP^{14}$ denote the
closed subscheme parametrizing pairs $(\Sigma,Q)$ where $Q\subset
\PP^4$ is a quadric hypersurface and $\Sigma\subset Q$.  The fiber of
the projection $I~\rightarrow~ \text{Hilb}_{P(t)}(\PP^4)$ over a
point $[\Sigma]$ is the projective space corresponding to
$H^0(\PP^4,\mathcal{I}_\Sigma(2))$, where $\mathcal{I}_\Sigma$ is the
ideal sheaf of $\Sigma\subset \PP^4$.  

\ps

Let $\widetilde{\Sigma}\subset \text{Hilb}_{P(t)}(\PP^4)\times \PP^4$
denote the universal closed subscheme, and let $\mathcal{I}$ denote the
ideal sheaf of this closed subscheme.  Consider the coherent sheaf
\begin{equation}\mathcal{F} ~=~
\text{pr}_{1*}(\mathcal{I}\otimes\text{pr}_2^*\OO_{\PP^4}(2)).\end{equation} 
For each $[\Sigma]\in \text{Hilb}_{P(t)}(\PP^4)$ there is an
evaluation map $\mathcal{F}|_{[\Sigma]} ~\rightarrow~
H^0(\PP^4,\mathcal{I}_{\Sigma}(2))$.  We will show that
$H^{i>0}(\PP^4,\mathcal{I}_{\Sigma}(2)) ~=~ 0$.  Then it follows by cohomology
and base change \cite[theorem III.12.11]{H} that $\mathcal{F}$ is locally free
in a neighborhood of $T$ and that all the evaluation maps are
isomorphisms in a neighborhood of $T$.  Thus in a neighborhood of $T$,
$I~\rightarrow~ \text{Hilb}_{P(t)}(\PP^4)$ is just the projective bundle
associated to $\mathcal{F}$.  

\ps

To show that $H^{i>0}(\PP^4,\mathcal{I}_{\Sigma}(2)) ~=~ 0$, we will use the short
exact sequence of coherent sheaves:
\begin{equation}\begin{CD}
0 @>>> \mathcal{I}_{\Sigma}(2) @>>> \OO_{\PP^4}(2) @>>> \OO_{\Sigma}(2) @>>> 0.
\end{CD}\end{equation}
Applying the long exact sequence in cohomology to this short exact
sequence, we see that we need to prove two things:
\begin{enumerate}
\item $H^{i>0}(\Sigma,\OO_{\Sigma}(2)) ~=~ 0$
\item $H^0(\PP^4,\OO_{\PP^4}(2)) ~\rightarrow~
H^0(\Sigma,\OO_{\Sigma}(2))$ is surjective.
\end{enumerate}

\ps

To prove (1) and (2) we will use the short exact sequence:
\begin{equation}\begin{CD}
0 @>>> \OO_{\Sigma_1}(2)(-L) @>>> \OO_{\Sigma}(2) @>>>
\OO_{\Sigma_2}(2) @>>> 0.
\end{CD}\end{equation}
Of course $\OO_{\Sigma_1}(2)(-L)\cong \OO_{\PP^2}(1)$ and
$\OO_{\Sigma_2}(2) \cong \OO_{\PP^1\times\PP^1}(2,2)$.  So applying the
long exact sequence in cohomology, we conclude that 
\begin{equation}H^1(\Sigma,\OO_{\Sigma}(2)) ~=~ H^2(\Sigma,\OO_{\Sigma}(2)) ~=~ 0\end{equation}
i.e., we have established (1).

\ps

To see that (2) is true, observe first that the composite map 
\begin{equation}\begin{CD} H^0(\PP^4,\OO_{\PP^4}(2)) @>>>
H^0(\Sigma,\OO_{\Sigma}(2)) @>>>
H^0(\Sigma_2,\OO_{\Sigma_2}(2))\end{CD}\end{equation} 
is surjective, i.e., the linear system $|\OO_{\PP^1\times\PP^1}(2,2)|$
on a smooth quadric surface is just the restriction of the linear
system $|\OO_{\PP^3}(2)|$.  The kernel of the composite map is
the vector space of quadratic polynomials which vanish identically on
$\text{span}(\Sigma_2)$.  If $F$ is a linear polynomial defining
$\text{span}(\Sigma_2)$, this subspace is just the image of the
multiplication map 
\begin{equation}\begin{CD} H^0(\PP^4,\OO_{\PP^4}(1)) @> *F >>
H^0(\PP^4,\OO_{\PP^4}(2)).
\end{CD}\end{equation}
The restriction to $H^0(\Sigma_1,\OO_{\Sigma_1}(2)(-L)) ~=~
H^0(\PP^2,\OO_{\PP^2}(1))$ is the restriction
$H^0(\PP^4,\OO_{\PP^4}(1))  \rightarrow  H^0(\PP^2,\OO_{\PP^2}(1))$,
which is clearly surjective.  Thus we have established (2).

\ps

We conclude that near $T$, $\mathcal{F}$ is locally free and
$I~\rightarrow~ \text{Hilb}_{P(t)}(\PP^4)$ is just the projective bundle
associated to $\mathcal{F}$.  If we let $V'$ be the open subset of 
$\text{Hilb}_{P(t)}(\PP^4)$ where 
$H^{i>0}(\PP^4,\mathcal{I}_{\Sigma}(2)) = 0$, then we know that $V$ is 
contained in $V'$, and that over $V'$ the map $I\rightarrow
\text{Hilb}_{P(t)}(\PP^4)$ is smooth (and hence flat).  This means that if
we have any open subset $O$ of $I$ over $V'$, its image in $V'$, and hence in
$\text{Hilb}_{P(t)}(\PP^4)$, will be open.

\ps

Notice that by the Lefschetz hyperplane theorem, there is no pair
$(\Sigma,Q)\in I$ such that $Q\subset \PP^4$ is a rank 5 quadric
(i.e., a smooth quadric).  Denote by $W\subset I$ the open
subscheme parametrizing pairs $(\Sigma,Q)$ such that $Q$ is a rank 4
quadric.  Denote by $W_T\subset W$ the locally closed subset such that
$\Sigma\in T$ and the singular point of $Q$ is a smooth point of
$\Sigma_1$. 

\ps

\begin{Claim}~\label{claim-1} The map $W_T~\rightarrow~ T$ is
surjective.
\end{Claim}
  
If $[\Sigma]$ is any point in $T$, and $p$ any point of $\Sigma_1$ not on
$\Sigma_2$, then letting $Q$ be the cone over $\Sigma_2$ with vertex $p$
provides a point of $W_T$ over $[\Sigma]$, which establishes the claim.
Now define $O\subset W$ to be the open
subset parametrizing pairs 
$(\Sigma,Q)$ such that the singular point of $Q$ is a smooth point of
$\Sigma$. 

\ps

\begin{Claim}~\label{claim-2}  $O$ is an irreducible open 
neighborhood of $W_T\subset W$ whose points $(\Sigma,Q)$ are exactly the 
points of $W_T$ and the pairs with $[\Sigma]\in U$.  
\end{Claim}

As part of the proof of the claim, we will see that there {\em are} points of
$O$ with $[\Sigma]\in U$. Since $U$ is a
homogeneous space for $\text{PGL}(5)$, and since $\text{PGL}(5)$ acts on the
rank four quadrics $Q$ as well, this means that the image of $O$ is exactly
$V=T\cup U$.  This will show that $V$ is open in
$\text{Hilb}_{P(t)}(\PP^4)$, since $O$ is open in $I$ over $V'$, and also
that $V$ is irreducible, since $O$ is. Finally, we know that $V$ is a smooth
subset of $\text{Hilb}_{P(t)}(\PP^4)$ since $U$ is smooth, and the Hilbert
scheme is smooth along $T$ (by lemma \ref{lem-specsm}).  Thus, the only step
left in proving lemma \ref{lem-smscr} is to establish claim \ref{claim-2}
above.

\ps

So we are reduced to studying the open neighborhood $O$.  
If we let $\tilde{Q}~\rightarrow~ Q$ denote the blow-up of $Q$ at $p$, and
if we let $\widetilde{\Sigma}\subset \tilde{Q}$ denote the proper
transform of $\Sigma$, then this open subset is also the parameter
space for pairs $(\widetilde{\Sigma},\tilde{Q})$.  

\ps

We will describe the
threefold $\tilde{Q}$.  Projection from $p$ defines a morphism
$\tilde{Q}~\rightarrow~ \PP^3$ whose image is a smooth quadric surface
$S\subset \PP^3$.  Identifying $S$ with $\PP^1\times \PP^1$, the projection
$\pi:\tilde{Q}~\rightarrow~ S$ is simply the $\PP^1$-bundle associated to
the vector bundle $G~=~\OO_{\PP^1\times\PP^1} \oplus \OO_{\PP^1\times
\PP^1}(1,1)$.  The exceptional divisor $E$ of $f:\tilde{Q}~\rightarrow~ Q$
is a section of $\pi$.  Identifying $E$ with $\PP^1\times \PP^1$, the
normal bundle of $E$ in $\tilde{Q}$ is identified with
$\OO_{\PP^1\times\PP^1}(-1,-1)$.  Let $F_1,F_2$ denote the divisor
classes of $\pi^*\text{pr}_1^*\OO_{\PP^1}(1), \pi^*\text{pr}_2^*\OO_{\PP^1}(1)$.
Then $\text{Pic}(\tilde{Q}) ~=~ \ZZ\{E,F_1,F_2\}$,  and for any $(\Sigma,Q)\in
O$, the proper transform $\widetilde{\Sigma}$ is a Cartier divisor with
divisor class $E+2F_1+F_2$ (up to permuting $F_1$ and $F_2$).  

\ps

Notice that since $p\in \Sigma$ is a smooth point, the intersection
$\widetilde{\Sigma}\cap E$ is a $(-1)$-curve in $\widetilde{\Sigma}$ along 
which $\widetilde{\Sigma}$ is smooth.
Conversely, suppose that $\Gamma\in |E+2F_1+F_2|$ is a surface
such that $\Gamma\cap E$ is a curve along which $\Gamma$ is smooth
(actually $\Gamma$ is automatically smooth along $\Gamma\cap E$ if
$\Gamma\cap E$ is a curve, but we won't need this fact).  We will show that
either $\Gamma$ is smooth or else 
$\Gamma$ is reducible, $\Gamma~=~\Gamma_1\cup \Gamma_2$ where $\Gamma_1$
is a smooth, connected divisor in the class of $F_1$,
$\Gamma_2$ is a smooth section of $\pi$ in the class of
$E+F_1+F_2$, and $\Gamma_1\cap \Gamma_2$ is transverse and maps to a
line in $Q$.  Then it follows that $f(\Gamma)$ is either a cubic
scroll or else the union $f(\Gamma_1)\cup f(\Gamma_2)$ of a 2-plane
and a smooth quadric surface along a line, and $p\in f(\Gamma_1)$ is a
smooth point.  

\ps

If $\Gamma$ is smooth, it is clear that $f(\Gamma)$ is a cubic scroll
(it is a smooth connected surface with Hilbert polynomial $P(t)$).
Therefore suppose that $\Gamma$ is singular at some point $q$.  We
know that 
$q\not\in E$.  

\ps

Suppose we pick a line $L$ in the quadric surface $S$, in the ruling 
corresponding to $F_1$.  If we restrict the $\PP^1$ bundle $\tilde{Q}$ over
$S$ to
$L$, the resulting surface is a Hirzebruch surface $\FF_1$ over $L$.  
A divisor
in the class $F_2$ on $\tilde{Q}$ restricts to the class of a fibre $F$ on 
$\FF_1$,
the exceptional divisor $E$ restricts to the unique $(-1)$-curve $D$, and a
divisor in the class $F_1$ restricts to the trivial class on $\FF_1$.

\ps

Now let $L_q$ be the particular line of ruling on $S$ containing $\pi(q)$, and
$\FF_1$ the surface over $L_q$.  If $\Gamma$ doesn't contain this $\FF_1$, then
the intersection $\Gamma\cap\FF_1$ is a curve on $\FF_1$ in the class
$|D+F|$, with a singular point at $q$, which is not on $D$.  This is a 
contradiction since the only singularities in the linear system $|D+F|$ occur
along $D$. (In the model of $\FF_1$ as the blowup of $\PP^2$ at a point, this
linear series is the pullback of the lines.)

\ps

Therefore the existence of a singular point $q \in\Gamma$ means that $\Gamma$
is reducible, and can be written $\Gamma_1\cup\Gamma_2$, with $\Gamma_1$ in 
the class $F_1$, and $\Gamma_2$ in the class $E+F_1+F_2$.  
Now since $N_{E/\tilde{Q}}\cong \OO_{\PP^1\times \PP^1}(-1,-1)$, we see that
$E\cap\Gamma_2$ is in the linear series $|\OO_{\PP^1\times\PP^1}|$, which means
that $\Gamma_2$ and $E$ are disjoint (we know that $\Gamma_2$ doesn't contain
$E$ as a component since $\Gamma$ intersects $E$ in a curve).
This shows both that the point $p$ lies on $f(\Gamma_1)$, and that if $\Gamma_2$
were to have a singular point $q'$, this point would not lie on $E$.

\ps

If $\Gamma_2$ were
to have 
any singular points, then the same argument as above would show that 
$\Gamma_2$ would be 
reducible, with one piece in the class of $F_1$ and one piece in the class of
$E+F_2$.  However, every element of $|E+F_2|$ contains $E$ as a component, which
is again a contradiction.
We conclude that $\Gamma_2$ is smooth.

\ps

The above arguments show that 
$f(\Gamma)$ is the union of a 2-plane $f(\Gamma_1)$ and a smooth
quadric surface $f(\Gamma_2)$ meeting along a line, and that $p$ lies on the
2-plane.  

\ps

We now know that every point in $O$ is either in $W_T$ or of
the form $(\Sigma,Q)$ with $\Sigma$ a cubic scroll.  Notice also that
$O$ fibers over the homogeneous space of rank 4 scrolls $Q$ and the
fiber over a point $[Q]$ is an open subset of the linear system
$|E+2F_1+F_2|$ on $\tilde{Q}$.  In particular the fibers are
irreducible, so $O$ is irreducible.   This finishes the proof of
claim \ref{claim-2}, and hence of
lemma \ref{lem-smscr}.

\end{proof}

\subsection{Additional constructions}~\label{subsec-addcon}

\ps

We prove several additional constructions of cubic scrolls which will
be needed.

\ps

Recall that our fourth description of a cubic scroll was an embedding
$f:\Sigma ~\rightarrow~ \PP^4$ where $\Sigma$ is the
Hirzebruch surface $\mathbb{F}_1$, and $f^*\OO(1)\sim \OO_{\Sigma}(1) ~=~
\OO_{\PP(E)}(D+2F)$.  The fact that the map is an embedding is equivalent to
asking that
the map $f$ be given by the
complete linear series of $\OO_{\Sigma}(D+2F)$.  In the next 
sections it
will be useful to weaken this condition.

\ps

\begin{defn} A \emph{cubic scroll} in $\PP^n$ is a finite morphism
$f:\Sigma~\rightarrow~ \PP^n$ where $\Sigma$ is isomorphic to the
Hirzebruch surface $\mathbb{F}_1$ and such that $f^*\OO_{\PP^n}(1)$ is
isomorphic to $\OO_{\Sigma}(D+2F)$. 
\end{defn}

\ps

We wish to look at various types of curves $C$ in $\PP^4$, and find conditions
for them to be enveloped by or contained in a cubic scroll $\Sigma$.  We always
start by looking at the class of the curve on $\Sigma$, look at its 
behavior with respect to the ruling and the directrix, and then seek to 
reconstruct the scroll out of this type of data.  When talking about the 
``degree'' of a curve $C$ on $\Sigma$, we always mean with respect to the line 
bundle $\OO_{\Sigma}(D+2F)$, which will be used to map $\Sigma$ into $\PP^4$.

\ps

\begin{lem}\label{lem-rib1}  Suppose $L\subset \PP^4$ is a line and
$T\subset T_{\PP^4}|_L$ is a sub-line bundle such that $T~\cong~
\OO_L(-1)$.  Then there is a unique scroll $f:\Sigma \rightarrow \PP^4$ with
$f(D)=L$ $(D$ the directrix of $\Sigma)$ and such that the differential map
$df_{*}: T_{\Sigma}\rightarrow f^{*}T_{\PP^4}$ takes the vertical tangent
bundle $T_{\Sigma/D} \subset T_{\Sigma}|_D$ to the pullback $f|_D^{*}T 
\subset f^{*}T_{\PP^4}|_D$ on $D$.
\end{lem}

\begin{proof}  We have the restriction to $L$ of the Euler sequence
for $\PP^4$:
\begin{equation}\begin{CD}
0 @>>> \OO_L @>>> \OO_L(1)^5 @>>> T_{\PP^4}|_L @>>> 0
\end{CD}\end{equation}
Define $F\subset \OO_L(1)^5$ to be the preimage of $T\subset
T_{\PP^4}|_L$.  We have an exact sequence:
\begin{equation}\begin{CD}
0 @>>> \OO_L @>>> F @>>> T\cong \OO_L(-1) @>>> 0
\end{CD}
\end{equation}
As $\text{Ext}^1_{\OO_L}(\OO_L(-1),\OO_L) ~=~ H^1(L,\OO_L(1)) ~=~ 0$, we
see that $F\cong \OO_L \oplus \OO_L(-1)$.  Therefore $E:~=~ F(-1)$ is
isomorphic to $\OO_L(-1)\oplus \OO_L(-2)$.  By construction $E$ is 
a subbundle $E~\hookrightarrow~ \OO_L^5$.  Defining $\Sigma~=~\PP(E)$, we have 
a map $\Sigma \rightarrow \PP(\OO_{L}^5)\cong \PP^1 \times \PP^4$.
Projecting onto the second factor, we get an induced map
$f:\Sigma~\rightarrow~ \PP^4$.

\ps

The directrix $D$ is the section of $\Sigma$ corresponding to
$\OO_L(-1)\subset E$.  The composite map $\OO_L(-1)\hookrightarrow 
E\hookrightarrow \OO_L^5$ is simply obtained from the first map of the 
Euler sequence by twisting by $\OO_L(-1)$,  and by construction of the 
Euler sequence, this is the same as the tautological map
$\OO_L(-1)~\rightarrow~ \OO_L^5$ induced by the inclusion $L\subset
\PP^4$.  Therefore $f(D)$ is just our original embedding of
$L$ in $\PP^4$.  Moreover, the vertical tangent bundle of $\Sigma$ on
$D$ is identified with
\begin{equation}
\text{Hom}(\OO_L(-1),E/\OO_L(-1))\subset
\text{Hom}(\OO_L(-1),\OO_L^5/\OO_L(-1))~=~ f^{*}T_{\PP^4}|_L .
\end{equation}
By construction of $E$ as the preimage of $T$, this is precisely 
$F/\OO_L ~=~ f^{*}T\subset f^{*}T_{\PP^4}|_D$.  
As well, $f^{*}\OO_{\PP^4}(1)=\OO_E(1)=\OO_{\Sigma}(D+2F)$, and 
therefore $f:\Sigma~\rightarrow~ \PP^4$ is the necessary scroll.

\ps

To see that $\Sigma$ is unique, simply observe that the lines through
the points of $L$ are determined by the direction of $T$ in
$T_{\PP^4}$.  Since the scroll is the union of its lines through $L$,
we see that $T$ uniquely determines the scroll.
\end{proof}

\ps

\begin{lem}\label{lem-rib2}  Let $C\subset \PP^4$ be a smooth conic
curve and let $T\subset T_{\PP^4}|_C$ be a sub-line bundle isomorphic
to $\OO_C(1)$ $($a degree $1$ line bundle on $C$, not
$\OO_{\PP^4}(1)|_C)$.  Then there is a unique scroll $f:\Sigma~\rightarrow~
\PP^4$ and a factorization $i:C\rightarrow \Sigma$ of $C\rightarrow \PP^4$ 
such that the differential $df:T_{\Sigma}\rightarrow f^{*}T_{\PP^4}$
maps the vertical tangent bundle $i^{*}T_{\Sigma/C}$
to $i^{*}f|_C^{*}T$ on $C$.
\end{lem}

\begin{proof}
As in the proof of lemma~\ref{lem-rib1}, define $F$ to be the
subbundle of $\OO_C(2)^5$ which is the preimage of $T\subset
T_{\PP^4}|_L$.  We have a short exact sequence:
\begin{equation}\begin{CD}\label{eqn-1000}
0 @>>> \OO_C @>>> F @>>> T\cong \OO_C(1) @>>> 0
\end{CD}\end{equation}
Since $H^1(C,\OO_C(-1))~=~0$, we have that $F\cong \OO_C\oplus
\OO_C(1)$.  Therefore $E:~=~ F(-2) \subset \OO_C^5$ is isomorphic to
$\OO_C(-2)\oplus \OO_C(-1)$.  Define $\Sigma ~=~ \PP(E)$.  The injective
map $E~\rightarrow~ \OO_C^5$ induces a morphism $f:\Sigma ~\rightarrow~
\PP^4$.  Define $i:C~\rightarrow~ \Sigma$ to be the section
associated to the twist by $\OO_C(-2)$ of the injection from
equation~(\ref{eqn-1000}): $\OO_C~\rightarrow~ F$.  

\ps

The composite map
$\OO_C(-2)~\rightarrow~ E ~\rightarrow~ \OO_C^5$ is just the twist of the
map in the Euler sequence $\OO_C~\rightarrow~ \OO_C(2)^5$.  By
construction of the Euler sequence, this is the map
$\OO_C(-2)~\rightarrow~ \OO_C^5$ corresponding to
$\OO_{\PP^4}(-1)~\rightarrow~ \OO_{\PP^4}^5$ induced by the inclusion
$C~\rightarrow~ \PP^4$.  Therefore $f(i(C))$ is just our original
embedding of $C$ in $\PP^4$.  Finally, notice that the restriction to
$i(C)$ of the vertical tangent bundle is simply 
\begin{equation}\begin{CD}
\text{Hom}(\OO_C(-2),E/\OO_C(-2)) @>>>
\text{Hom}(\OO_C(-2),\OO_C^5/\OO_C(-2)) ~=~ f^{*}T_{\PP^4}|_C
\end{CD}\end{equation}
By construction of $E$, this is precisely $F/\OO_C ~=~ f^{*}T\subset
f^{*}T_{\PP^4}|_C$.  Finally, 
$f^{*}\OO_{\PP^4}(1) = \OO_E(1) = \OO_{\Sigma}(D+2F)$, 
and therefore $f:\Sigma~\rightarrow~ \PP^4$ is the necessary scroll.

\ps

To see that $\Sigma$ is uniquely determined, observe that $T\subset
T_{\PP^4}|_C$ determines the lines in $\Sigma$ which pass through
$C$.  Since $\Sigma$ is the union of the lines which pass through $C$, 
this shows that $\Sigma$ is unique.
\end{proof}

\ps

\subsection{Cubic Scrolls and Quartic Rational
Curves}~\label{subsec-scr40}

\ps

Recall that $\text{Pic}(\Sigma) = \ZZ\{D,F\}$ where $D$ is the
directrix and $F$ is the class of a line of ruling.  The intersection
product is given by $D^2 = -1$, $D.F = 1$, $F^2 =0$.  The canonical
class is given by $K_\Sigma = -2D-3F$.  Our definition of a cubic scroll
is that the finite map $f:\Sigma\rightarrow \PP^4$ should come from the line
bundle $f^{*}\OO_{\PP^n}(1)=\OO_{\Sigma}(D+2F)$. 

\ps

The linear system $|F|$ is nef because it is the
pullback of $\OO_{\PP^1}(1)$ under the projection
$\pi:\Sigma\rightarrow \PP^1$.  Similarly, $|D+F|$ is nef because 
it is the pullback of $\OO_{\PP^2}(1)$ in the
realization of $\Sigma$ as $\PP^2$ blown up at a point.
Thus for any effective curve class $aD+bF$ we have
the two inequalities $a=(aD+bF).F \geq 0$, $b=(aD+bF).(D+F) \geq 0$.  

\ps

Suppose that $C\subset \Sigma$ is an effective
divisor of degree $4$ and arithmetic genus $0$.  By the adjunction
formula 
\begin{equation}
K_\Sigma.[C] + [C].[C] = 2p_a - 2 = -2.
\end{equation}
So if $[C]=aD+bF$, then we have the conditions
\begin{equation}
a\geq 0,\ b\geq 0,\ a+b = 4,\ a^2 -2ab +a+2b = 2.
\end{equation}
It is easy to check that there are precisely two solutions $[C]=2D+2F$, and
$[C]=D+3F$.  We will see that both possibilities occur and describe
some constructions related to each possibility.

\ps

We start with the case $[C] = 2D + 2F$.

\ps

\begin{lem}\label{lem-g1240}  Let $C\subset \PP^4$ be a smooth quartic 
rational curve and let $V\subset |\OO_C(2)|$ be a pencil of degree
$2$-divisors on $C$ without basepoints.  There exists a unique cubic scroll
$f:\Sigma~\rightarrow~ \PP^4$ 
and a factorization $i:C~\rightarrow~ \Sigma$ of $C~\rightarrow~ \PP^4$
such that $[i(C)] = 2D+2F$ and such that the pencil of degree $2$ divisors
$\pi^{-1}(t)\cap C$ $($for $t\in\PP^1)$ is the pencil $V$. 
\end{lem}

\begin{proof}
Let $g:C~\rightarrow~ \PP^1$ be a degree $2$-morphism defining $V$.
Define $E^\vee :~=~ g_* (\OO_{\PP^4}(1)|_C)$.  Since $g$ is finite and
flat of degree $2$, $E^\vee$ is locally free of rank $2$.
Since $g^{*}\OO_{\PP^1}(1)\cong \OO_{C}(2)$, the projection formula shows
that $E^\vee\cong \OO_{\PP^1}(2)\otimes g_*\OO_C$.  But
$g:C~\rightarrow~ \PP^1$ is a cyclic cover of degree $2$ branched over a 
divisor of degree $2$.  The theory of cyclic covers~\cite[definition
2.50]{KM} shows that $g_*\OO_C$ decomposes as a sum of
$\ZZ/2\ZZ$-eigensheaves: $g_*\OO_C \cong
\OO_{\PP^1} \oplus \OO_{\PP^1}(-1)$.  Thus $E\cong
\OO_{\PP^1}(-2)\oplus \OO_{\PP^1}(-1)$.  

\ps

On $\PP^4$ we have the surjection of vector bundles $\OO_{\PP^4}^5\rightarrow
\OO_{\PP^4}(1)$ given by global sections of $\OO_{\PP^4}(1)$.  Restricting
this to $C$ gives $\OO_{C}^5~\rightarrow~ \OO_{\PP^4}(1)|_C$, 
which by adjunction induces a map 
$\OO_{\PP^1}^5~\rightarrow~ g_{*}\OO_{C}(2) = E^\vee$.  
Since $C~\rightarrow~ \PP^4$ is an
embedding, for each pair of points $\{p,q\}\subset C$ (possibly
infinitely near), we have that $\OO_C^5~\rightarrow~
\OO_{\PP^4}(1)|_{\{p,q\}}$ is surjective.  In particular, taking
$\{p,q\} ~=~ g^{-1}(t)$ for $t\in\PP^1$, we conclude that
$\OO_{\PP^1}^5~\rightarrow~ E^\vee|_t$ is surjective.  Define $\Sigma ~=~
\PP(E)$, then the surjective map $\OO_{\PP^1}^5~\rightarrow~ E^\vee$
induces a morphism $f:\Sigma~\rightarrow~ \PP^4$.  We have $E\cong
\OO_{\PP^1}(-2)\oplus \OO_{\PP^1}(-1)$ so that $\Sigma \cong \FF_1$,
and $f^*\OO_{\PP^4}(1) ~=~ \OO_{E}(1) =  \OO_{\Sigma}(D+2F)$.

\ps

By adjunction we have a map of sheaves $g^*E^\vee ~\rightarrow~
\OO_{\PP^4}(1)|_C$.  This map is surjective since $g$ is finite.
Moreover $g^* E^\vee ~\rightarrow~ \OO_{\PP^4}(1)|_C$ even separates
points: for points in distinct fibers this is clear.  For points
$\{p,q\}~=~g^{-1}(t)$ it follows because $E^{\vee}|_t$ is precisely
$\OO_{\PP^4}(1)|_{\{p,q\}}$.  So the induced morphism $i:C~\rightarrow~
\Sigma$ is even an embedding.  Moreover, the pullback
$i^*\OO_{\Sigma}^5~\rightarrow~ i^*\OO_{E}(1)$ is precisely our
original map $\OO_C^5~\rightarrow~ \OO_{\PP^4}(1)|_C$ so that $f\circ
i:C~\rightarrow~ \PP^4$ is our original embedding of $C$ in $\PP^4$.

\ps

By construction of
$E$, we have $i^*|F|~=~V$,  and $g^*\OO_{\PP^1}(-1)~=~\OO_C(-2)$ so we see 
$i^*\OO_{\Sigma}(D)\cong \OO_C(2)$.  Thus we have $i(C)\sim 2D+2F$.
Therefore $f:\Sigma~\rightarrow~ \PP^4$, $i:C~\rightarrow~ \Sigma$ are the
necessary maps.

\ps

The map $f$ is only an embedding if $C$ is nondegenerate,
otherwise $f$ is the normalization map for its image, which is a singular 
cubic surface.

\ps

To see that this is unique, notice that the lines $f(\pi^{-1}(t))$ are simply
the lines obtained by taking the joins of the degree $2$ divisors on
$C$ which lie in in $V$.  Since $f(\Sigma)$ is the union of this system
of lines, this proves that $f(\Sigma)$ is uniquely determined.  But
$f:\Sigma~\rightarrow~ f(\Sigma)$ is simply the normalization map so that
$f$ is also uniquely determined.
\end{proof}

\textbf{Remark}
While we are at it, let's mention a specialization of the construction
above, namely what happens when $V$ is {\em not} basepoint free.
Then $V=p+|\OO_C(1)|$, where $p\in C$ is some basepoint.  Consider
the projection morphism $f:\PP^4\dottedrightarrow \PP^3$ obtained by
projection from $p$ (this is a rational map undefined at $p$).  The
image of $C$ is a rational cubic curve $B$ (possibly a singular plane
cubic).  Consider the cone $S'$ in $\PP^4$ over $B$ with vertex
$p$.  This 
surface contains $C$.  If we blowup $\PP^4$ at $p$, then the proper
transform of $S'$ in $\widetilde{\PP^4}$ is a surface whose
normalization $S$ is a 
Hirzebruch surface $\FF_3$ (normalization is only necessary if
$B$ is a plane curve).  The pullback of the exceptional divisor of 
$\widetilde{\PP^4}$ plays the role of the directrix $D$ of $S$.  The inclusion
$C\subset S'$ induces a factorization $i:C\rightarrow S$ of
$C\rightarrow \PP^4$, with $[i(C)] = D + 4F$.  The intersection of $D$
and $i(C)$ is precisely the point $p$.  Finally, the linear system $i^*|F|$
is exactly $|\OO_C(1)|$.  

\ps

Next we consider the case of a rational curve $C\subset \Sigma$ such
that $[C]=D+3F$.  

\ps

\begin{lem}~\label{lem-scr40}  
Let $C\subset \PP^4$ be a smooth quartic rational curve 
and let $L\subset \PP^4$ be a line such that $L\cap C=Z$ is a degree 2
divisor.  
Let $\phi:C\rightarrow L$ be an isomorphism such that $\phi(Z)=Z$ and
$\phi|_Z$ is the identity map.  Then there exists a unique triple
$(f,i,j)$ where $f:\Sigma\rightarrow \PP^4$ is a cubic scroll,
$i:C~\rightarrow~\Sigma$ and $j:L~\rightarrow~\Sigma$ are factorizations of
$C\rightarrow \PP^4$ and $L\rightarrow \PP^4$, and 
such that $j(L)=D$ is the directrix, $[i(C)]=D+3F$, and 
the lines of ruling induce the original isomorphism $\phi:C\rightarrow L$.
\end{lem}

\begin{proof}
Choose isomorphisms $g:\PP^1\rightarrow C$ and $h:\PP^1\rightarrow L$
such that $\phi\circ g = h$.  Consider the rank 2 vector bundle
\begin{equation}
G=g^*\OO_{\PP^4}(1) \oplus h^*\OO(1) \cong
\OO_{\PP^1}(4)\oplus\OO_{\PP^1}(1). 
\end{equation}
Let $Z=g^*Z = h^*Z$.  Since $g(Z)=h(Z)$ as subschemes of $\PP^4$, we
have an identification of $g^*\OO_{\PP^4}(1)|_Z$ with
$h^*\OO_{\PP^4}(1)|_Z$.  Define $E^\vee\subset G$ to be the subsheaf
of $E$ of sections $(s_C,s_L)$ such that $s_C|_Z=s_L|_Z$ under our
identification. 

\ps

Since the map $g^*\OO_{\PP^4}(1) \rightarrow g^{*}\OO_{\PP^4}(1)|_Z$ 
is surjective, we conclude that
$E^\vee~\cong~\OO_{\PP^1}(2)\oplus\OO_{\PP^1}(1)$.  Moreover the
linear series 
\begin{equation}
H^0(\PP^4,\OO_{\PP^4}(1))\rightarrow H^0(\PP^1,g^*\OO_{\PP^4}(1)\oplus
h^*\OO_{\PP^4}(1)) 
\end{equation}
clearly factors through $H^0(\PP^1,E^\vee)$.  The question arises
whether 
\begin{equation}
H^0(\PP^4,\OO_{\PP^4}(1))\otimes_\CC \OO_{\PP^1}\rightarrow
E^\vee = \OO_{\PP^1}(2)\oplus\OO_{\PP^1}(1)
\end{equation}
is surjective.  Certainly the corresponding maps to
$g^*\OO_{\PP^4}(1)$ and $h^*\OO_{\PP^4}(1)$ are surjective.  The
condition that $C\cap L = Z$ is precisely the condition that the image
of 
\begin{equation}
H^0(\PP^4,\OO_{\PP^4}(1))\otimes_\CC \OO_{\PP^1}\rightarrow
f^*\OO_{\PP^4}(1)\oplus g^*\OO_{\PP^4}(1)
\end{equation}
is $E^\vee$.  Thus the morphism is surjective.  Denoting $\Sigma=\PP(E)$,
we conclude that there is a well-defined morphism
$f:\Sigma\rightarrow \PP^4$ such that $f^*\OO_{\PP^4}(1)=\OO_{E}
(1)$ and the pullback map
\begin{equation}
H^0(\PP^4,\OO_{\PP^4}(1))\rightarrow H^0(\PP(E),\OO_{E}
(1))=H^0(\PP^1,E^\vee)
\end{equation}
is the map above.  

\ps

The composition of $E^\vee~\rightarrow~g^*\OO_{\PP^4}(1)\oplus
h^*\OO_{\PP^4}(1)$ with the two projections define two surjective maps
which yield sections $i:\PP^1~\rightarrow~ \Sigma$,
$j:\PP^1~\rightarrow~\Sigma$.  Clearly $\pi\circ i = g^{-1}$, and 
$\pi\circ j=h^{-1}$.
From this it follows that the isomorphism $C\cong L$ induced by the ruling
of $\Sigma$ is the same as the isomorphism $\phi$.  Thus $(f,i,j)$ is a 
triple as in the statement of the lemma.  

\ps

To see that this is unique, notice that the lines $f(\pi^{-1}(t))$ are simply
the lines obtained by $\text{span}(p,\phi(p))$.  Since $f(\Sigma)$ is
the union of this system of lines, we conclude that $f(\Sigma)$ is
unique.  But $f:\Sigma\rightarrow f(\Sigma)$ is just the normalization map
so that $f$ is also uniquely determined.
\end{proof}

\ps

\subsection{Cubic Scrolls and Quintic Elliptics}~\label{subsec-scr51}

\ps

Suppose that $E\subset \Sigma$ is an effective Cartier divisor
of degree $5$ and arithmetic genus $1$.
Writing $E~=~aD + bF$ we see
$(a,b)$ satisfies the relations $a,\, b \geq 0$, $a+b ~=~ 5$ and
$a(b-3)+b(a-2) -a(a-2) ~=~ 0$.  These relations give the unique integer solution
$E~=~2D+3F ~=~ -K$.  In particular, if $E$ is smooth then
$\pi:E~\rightarrow~ \PP^1$ is a finite morphism of degree 2, i.e., a
$g^1_2$ on $E$.  Thus a pair $(f:\Sigma~\rightarrow~ \PP^n, E\subset
\Sigma)$ of a cubic scroll and a quintic elliptic determines a pair
$(g:E~\rightarrow~ \PP^n, \pi:E~\rightarrow~ \PP^1)$ where $g:E~\rightarrow~
\PP^n$ is a quintic elliptic and $\pi:E~\rightarrow~ \PP^1$ is a degree
2 morphism.  

\ps

Suppose we start with a pair $(h:E~\rightarrow~ \PP^n$, $\pi:E~\rightarrow~
\PP^1)$ where $h:E~\rightarrow~ \PP^n$ is an embedding of a quintic
elliptic curve and $\pi:E~\rightarrow~ \PP^1$ is a degree 2 morphism.
Consider the rank 2 vector bundle $\pi_* h^* \OO_{\PP^n}(1)$.  

\ps

\begin{lem} Suppose $E$ is an elliptic curve and $\pi:E~\rightarrow~
\PP^1$ is a degree 2 morphism.  Suppose $L$ is an invertible sheaf on
$E$ of degree $d$.  Then we have

\begin{equation}\pi_* L \cong \lt\{ \begin{array}{ll}
 \OO_{\PP^1}(e) \oplus \OO_{\PP^1}(e-1) & d ~=~ 2e+1, \\
 \OO_{\PP^1}(e) \oplus \OO_{\PP^1}(e-2) & d~=~2e, L\cong
	\pi^*\OO_{\PP^1}(e), \\
 \OO_{\PP^1}(e-1) \oplus \OO_{\PP^1}(e-1) & d~=~2e, L\not\cong \pi^*
	\OO_{\PP^1}(e) 
\end{array} \rt.\end{equation}

\end{lem}

\begin{proof}
Using the projection formula, we see that the lemma for $L$ is
equivalent to the lemma for $L\otimes \pi^*\OO_{\PP^1}(m)$.  For
any $L$ there is an $m$ such that $L\otimes \pi^*\OO_{\PP^1}(m)$ has
degree 0 or degree 1.  Thus we are reduced to the two cases $d~=~0$ and
$d~=~1$.  Notice also that in all cases we have $\chi(\pi_* L) ~=~
\chi(L)$, so that by Riemann-Roch for $E$ and $\PP^1$ we have 
\begin{equation} 
\text{deg}(\pi_* L) + \text{rank}(\pi_* L) ~=~ \text{deg}(L) ~=~ d,
\end{equation}
i.e., $\text{deg}(\pi_* L) ~=~ d-2$.  By Grothendieck's lemma about
vector bundles on $\PP^1$ we know $\pi_* L ~=~ \OO_{\PP^1}(a)\oplus
\OO_{\PP^1}(d-2-a)$.  

\ps

Suppose now that $d~=~0$.  Then $\pi_*L ~=~ \OO_{\PP^1}(a)\oplus
\OO_{\PP^1}(-2-a).$  We also have $h^0(\pi_* L) = h^0(L)$.  If $L
\cong \OO_E$ then $h^0(L) ~=~ 1$ so that we have $\pi_* L ~=~
\OO_{\PP^1}\oplus \OO_{\PP^1}(-2)$.  If $L\not\cong \OO_E$ then
$h^0(L)~=~0$ so that we have $\pi_* L ~=~ \OO_{\PP^1}(-1)\oplus
\OO_{\PP^1}(-1).$  Thus the lemma is proved when $d~=~0$.

\ps

Next suppose that $d~=~1$.  Then $\pi_* L ~=~ \OO_{\PP^1}(a)\oplus
\OO_{\PP^1}(-1-a)$,  and $h^0(\pi_* L) ~=~ h^0(L) ~=~ 1$.  Thus we have
$\pi_* L ~=~ \OO_{\PP^1} \oplus \OO_{\PP^1}(-1).$  So the lemma is
proved when $d~=~1$.  Thus the lemma is proved in all cases.
\end{proof}

By the lemma we see that the vector bundle $G^{\vee}:~=~\pi_* h^*
\OO_{\PP^n}(1)$ is isomorphic to $\OO_{\PP^1}(1) \oplus 
\OO_{\PP^1}(2)$.  Associated to the linear series
$\OO_E^{n+1}~\rightarrow~ h^*\OO_{\PP^n}(1)$ defining the embedding $h$,
we have the adjoint map $\OO_{\PP^1}^{n+1}~\rightarrow~ 
\pi_{*}h^{*}\OO_{\PP^4}(1)= G^{\vee}$.  Since $h$ is an embedding, for each pair of points $\{p,q\}
\subset E$ (possibly infinitely near), we have that
$\OO_E^{n+1}~\rightarrow~ h^*\OO_{\PP^n}(1)|_{\{p,q\}}$ is surjective.
In particular taking $\{p,q\}~=~\pi^{-1}(t)$ for $t\in \PP^1$, we
conclude that $\OO_{\PP^1}^{n+1}~\rightarrow~ G^{\vee}|_t$ is surjective.  Thus
we have an induced morphism $\PP(G)~\rightarrow~ \PP^n$ which pulls
back $\OO_{\PP^n}(1)$ to $\OO_{G}(1)$.  Let us
denote $\Sigma :~=~ \PP(G)$ and let us denote the morphism by
$f:\Sigma~\rightarrow~ \PP^n$.  Abstractly $\Sigma$ is isomorphic to
$\mathbb{F}_1$ and $f:\Sigma~\rightarrow~ \PP^n$ is a cubic scroll.

\ps

The tautological map $\pi^*\pi_* h^*\OO_{\PP^n}(1)~\rightarrow~
h^*\OO_{\PP^n}(1)$ is clearly surjective.  Thus there is an induced
morphism $i:E~\rightarrow~ \Sigma$.  Chasing definitions, we see that
$h~=~f\circ i$.  So we conclude that given a pair $(h:E~\rightarrow~
\PP^n, \pi:E~\rightarrow~ \PP^1)$ as above, we obtain a pair
$(f:\Sigma~\rightarrow~ \PP^n, i:E~\rightarrow~ \Sigma)$.  Thus we
have proved the following:

\begin{lem}\label{lem-equiv1}
There is an equivalence between the collection of pairs
$(f:\Sigma~\rightarrow~ \PP^n, i:E~\rightarrow~ \Sigma)$ with
$f:\Sigma~\rightarrow~ \PP^n$ a cubic scroll and $f\circ i:E~\rightarrow~
\PP^n$ an embedded quintic elliptic curve and the collection of pairs
$(h:E~\rightarrow~ \PP^n, \pi:E~\rightarrow~ \PP^1)$ where $h:E~\rightarrow~
\PP^n$ is an embedded quintic elliptic curve and $\pi:E~\rightarrow~
\PP^1$ is a degree 2 morphism.  
\end{lem}

Stated more precisely, this gives an
isomorphism of the parameter schemes of such pairs, but we won't need
the result in this form.

\ps

\subsection{Cubic Scrolls and Quintic Rational
Curves}~\label{subsec-scr50}

\ps

If one carries out the analogous computations as at the beginning of
subsection~\ref{subsec-scr40} one sees that the only effective divisor
classes $aD+bF$ on a cubic scroll $\Sigma$ with degree $5$ and
arithmetic genus $0$ are $D+4F$ and $3D+2F$.  But the divisor class
$3D+2F$ cannot be the divisor of an irreducible curve because
$(3D+2F).D = -1$.  Thus if $C\subset \Sigma$ is an irreducible curve
of degree $5$ and arithmetic genus $0$, then $[C]=D+4F$.  

\ps

\begin{lem}~\label{lem-scr50}
Let $C\subset \PP^4$ be a smooth quintic rational curve
and let $L\subset \PP^4$ be a line such that $L\cap C$ is a degree 3
divisor $Z$.  Let $\phi:C\rightarrow L$ be an isomorphism such that
$\phi(Z)=Z$ and $\phi|_Z$ is the identity map.  Then there exists a
unique triple $(f,i,j)$ such that $f:\Sigma~\rightarrow~\PP^4$ is a
cubic scroll, $i:C~\rightarrow~\Sigma$, and $j:L~\rightarrow~\Sigma$ 
are factorizations of $C~\rightarrow~\PP^4$, and $L~\rightarrow~\PP^4$,
and such that $j(L)=D$ is the
directrix, $[i(C)]=D+4F$, and the lines of ruling induce the original 
isomorphism $\phi:C\rightarrow L$.
\end{lem}

\begin{proof}
The proof is very similar to the proof of lemma~\ref{lem-scr40}.  
\end{proof}

\ps

\section{Quartic Rational Curves}\label{sec-irr40} 

\ps

In this section we will prove that the space $\Hdg{4}{0}(X)$ of smooth 
quartic rational curves $C\subset X$ is irreducible of dimension $8$.
Recall from section~\ref{sec-def} that every irreducible component of
$\Hdg{4}{0}(X)$ has dimension at least $-K_X.C = 2\times 4 =8$.  First
we prove that the 
open subset $\mathcal{U}\subset \Hdg{4}{0}(X)$ parametrizing curves
$C$ with $\text{span}(C)=\PP^4$ is Zariski dense.  To prove this it
suffices to prove that the complement
$\mathcal{D}\subset\Hdg{4}{0}(X)$ has dimension at most $7$.
  
\ps

\begin{lem}~\label{lem-irr40}  Every irreducible component of the
closed subset $\mathcal{D}\subset 
\Hdg{4}{0}(X)$ parametrizing degenerate curves $C$
$($i.e., $\text{span}(C)\neq \PP^4)$ has dimension at most $7$.
\end{lem}

\begin{proof}
By Riemann-Roch we see that a smooth quartic rational curve $C\subset
\PP^3$ lies on at least one quadric surface $S$.  It cannot lie on two 
distinct quadric surfaces, for then it would have arithmetic genus $1$ 
which is a contradiction.  Thus to each point $[C]\in \mathcal{D}$,
there is an associated quadric surface $S\subset \text{span}(C)$.
Moreover the residual to $C\subset S\cap X$ is a pair of lines $L_1,
L_2$ (possibly a single non-reduced line).  Thus there is a morphism
$\mathcal{D} ~\rightarrow~ \HI{2}{F}$.  Since $F$ is a smooth surface,
$\Hi{2}(F)$ has dimension $4$.  

\ps

Now there are two types of behaviors depending on whether or not
$\text{span}(L_1,L_2)$ is a $2$-plane or a $3$-plane.  The set of
pairs $\{L_1,L_2\}$ such that $\text{span}(L_1,L_2)$ is a $2$-plane
corresponds to a point in the $3$-dimensional divisor $\mathcal{I}
\subset \Sym^2(F)$ parameterizing incident lines.
For each pair $\{L_1,L_2\}$ on this $3$-fold, there
is a $1$-parameter family of hyperplanes containing
$\text{span}(L_1,L_2)$.  For each such hyperplane, there is a $\PP^3$
of quadric surfaces $Q$ in this hyperplane which contain $L_1\cup
L_2$.  Thus the locus of all curves $[C]\in \mathcal{D}$ whose
associated pair $\{L_1,L_2\}$ lies in $\mathcal{I}$ has dimension at
most $3+1+3~=~7$.  

\ps

Next consider the case that $\text{span}(L_1,L_2)$ is a $3$-plane.
Then every quadric surface containing these lines lies in this
$3$-plane.  The set of quadric surfaces in this $3$-plane which
contain $L_1$ and $L_2$ is itself a $\PP^3$.  Thus the set of curves
$[C]\in \mathcal{D}$ whose associated pair $\{L_1,L_2\}$ spans a
$3$-plane has dimension at most $4+3~=~7$.  So the lemma is proved.
\end{proof}

Recall from lemma~\ref{lem-enum} that given any smooth quartic
rational curve $C\subset X$, the subscheme $A_C\subset
\text{Grass}(2,5)$ parametrizing the 2-secant lines to $C$ is either
1-dimensional or else is 0-dimensional of length 16.  In either case
we conclude that there exists a 0-dimensional, length 2 subscheme
$Z\subset A_C$ (in fact many such subschemes).
Suppose given a 0-dimensional, length 2 subscheme $Z\subset A_C$.
Then $Z$ either consists of two reduced points $[L_1],[L_2]\in A_C$ or 
else $Z$ corresponds to a non-reduced point of $A_C$.  
In the case that $Z=\{[L_1],[L_2]\}$, there are again two behaviors
depending on whether $Z$ is \emph{planar}, i.e., $\text{span}(L_1,L_2)$ 
is a 2-plane, or whether $Z$ is \emph{non-planar},
i.e., $\text{span}(L_1,L_2)$ is a 3-plane.  In the case that $Z$ is
planar, notice that we have the \emph{distinguished point} $p\in X$
corresponding to the intersection of $L_1$ and $L_2$.  
In order to explain the
analogues of planar and non-planar in the case that $Z$ is
non-reduced, we make a brief digression on ribbons.

\ps

A {\em ribbon} (for our purposes) is a degree two subscheme $R$ of $\PP^4$
supported along a line $L$, such that the ideal sheaf $\mathcal{I}_L$ of 
$L$ in $R$ satisfies $\mathcal{I}_L^2=0$, and such that the conormal sheaf
$\mathcal{I}_L/\mathcal{I}_L^2 = \mathcal{I}_L$ is a line bundle
on $L$.  A ribbon is therefore a kind of doubled line in $\PP^4$, such that
at each point $p$ of $L$ the doubling occurs in a specified direction (given
by the normal bundle $N_{L/R}=\mathcal{I}_L^{\vee}$), and such that this 
doubling direction varies reasonably along the line.
Starting with a fixed line $L$ in $\PP^4$, and a sub-line-bundle
$N\subset N_{L/\PP^4}$, there is a unique ribbon $R$ supported on $L$ with
$N_{L/R}=N$.

\ps

A non-reduced
0-dimensional subscheme $Z$ of $\text{Grass}(2,5)$ of length 2 
determines a ribbon $R$ in $\PP^4$.  The line $L$ of the ribbon is given by 
the point of support of $Z$ in $\text{Grass}(2,5)$, while the tangent 
direction of $Z$ corresponds to a global section of 
$N_{L/\PP^4}$, and there is a unique sub-line-bundle 
$N$ of $N_{L/\PP^4}$ containing this section, which gives the ribbon.

\ps
There are two possibilities.  

\ps

First of all we could have
$N_{L/R}~\cong~\OO_L(1)$.  In this case we say that $R$ is
\emph{planar ribbon} since there is a unique 2-plane $P\subset \PP^4$
such that $R\subset P$ -- in fact $P$ is the unique 2-plane such that
$N_{L/P} = N_{L/R}$ as subbundles of $N_{L/\PP^4}$.  Notice that in
this case the global section $\OO_L\rightarrow N_{L/R}$ is not
determined by the ribbon $R$ -- in fact the data of this section is
equivalent to a point $p\in L$ such that the length 2 scheme $Z$ is simply 
the tangent direction at $[L]$ to the pencil of lines in $P$ which
pass through $p$.  We refer to the point $p\in L$ as the
\emph{distinguished point} of $L$ determined by $Z$.

\ps

The second possibility for the ribbon is that $N_{L/R}~\cong~\OO_L$.
First of all notice that in this case $Z$ is uniquely determined by the ribbon. 
Second, given any subbundle $\OO_L \cong N \subset N_{L/\PP^4} \cong\OO_L(1)^3$,
there is a unique 
3-plane $H\subset \PP^4$ such that the map $\OO_L~\rightarrow~
N_{L/\PP^4}$ factors through $N_{L/H}~\subset~N_{L/\PP^4}$.  Moreover
in $H$ there is a $\PP^3$ of quadric surfaces $Q\subset H$ such that
$R\subset Q$.  The general surface $Q$ in this $\PP^3$  will be smooth 
and we will have $N_{L/R}=N_{L/Q}$ as subbundles of $N_{L/H}$.  We
will call a ribbon of this type a \emph{non-planar ribbon}.

\ps

Define $I=I_{4,0}\subset \mathcal{U}\times \HI{2}F$ to be the incidence
correspondence of pairs $([C],[Z])$ such that $Z\subset A_C$ is a
0-dimensional length 2 subscheme.  The idea of the proof of
irreducibility of $\Hdg{4}{0}(X)$ is to consider for such a pair
$([C],[Z])$ a certain cubic surface $\Sigma$ which contains the curve
which is the union of $C$ and the scheme parametrized by $Z$.  The
residual of this curve in $\Sigma\cap X$ will be a cubic curve, and
for general $([C],[Z])$ this will be a smooth cubic rational curve.
Moreover, if we associate to this cubic rational curve its image in
$\Theta$ under the Abel-Jacobi map, then we obtain a rational
transformation $I\rightarrow \Theta \times \HI{2}{F}$ as the
product of this map with projection $I\subset \Hdg{4}{0}(X)\times
\HI{2}{F} \rightarrow \HI{2}{F}$.  The main fact is that this rational 
transformation is birational.  

\

In order to prove the claims made in the last paragraph, we must first 
dispense with some degenerate possibilities.  Let $I_P\subset I$
denote the closed subset parametrizing pairs $([C],[Z])$ such that $Z$ 
is planar.  

\begin{lem}\label{lem-planar}  Every irreducible component of $I_P$
has dimension at most 7.
\end{lem}

\begin{proof}
In the reduced case
$Z~=~\{[L_1],[L_2]\}$, we have that $C\cap(L_1\cup L_2)$ is a
0-dimensional subscheme of length 4 unless the distinguished point
$p\in C$, in which case $C\cap(L_1\cup L_2)$ has length 3.  Similarly
in the case that $Z$ is non-reduced and gives rise to a planar ribbon
$R$, we have $C\cap R$ is length 4 unless the distinguished point
$p$ is on $C$.  But since $\text{span}(C)=\PP^4$, there is no 2-plane $P$
such that $P\cap C$ has length 4; if such a 2-plane exists, then for
any point $q\in C-P\cap C$ we have the hyperplane $H=\text{span}(P,q)$ 
intersects $C$ in a scheme of length 5 which contradicts B\'ezout's
theorem since $C\not\subset H$.  Therefore we conclude that if
$Z\subset A_C$ is planar, then the distinguished point $p$ lies on
$C$.  

\ps

Now define $S$ to be the cone over $C$ with vertex $p$.  The
projection of $\PP^4$ away from $p$ (to $\PP^4/p~\cong~\PP^3$) maps
$C$ birationally to a smooth cubic rational curve $C'$, and $S$
is simply the cone over this cubic curve.  Moreover, $S$ contains 
the curve $E$ which is the union of $C$ and the degree 2 subscheme
parametrized by $Z$ (either $L_1\cup L_2$ or else the ribbon $R$
determined by $Z$).  By the Lefschetz hyperplane theorem, $X$ contains 
no cubic surfaces other than the (degenerate) hyperplane sections of
$X$.  Now $S$ is non-degenerate since it contains $C$ and $C$ is
non-degenerate.  Therefore $S$ is not contained in $X$.  So
$S\cap X$ is a divisor on $S$ of degree $\text{deg}(S)
\times \text{deg}(X) = 3\times 3 = 9$.  But $E$ has degree $6$.  Thus
the residual curve $D$ to $E$ is a curve of degree $3$.  The only
curves of degree $3$ on $S$ are hyperplane sections.  There are
two possible cases depending on whether or not $p\in D$.

\ps

Suppose that $p\in D$.  In this case $D$ is a union of three lines in
$S$ through $p$ (or some degeneration thereof).  Let $H\subset
\PP^4$ be the tangent hyperplane to $X$ at $p$.  Then every line
$L\subset X$ containing $p$ is contained in $H$.  In particular, the
residual subscheme to $C$ in $S\cap X$ is contained in $H$.  But
an easy divisor class calculation on the blowup of $S$ at $p$
shows that the residual to $C$ intersects $C$ in a divisor of degree
$5$ (not counting $p$ where the residual isn't well-defined).  So
$C\cap H$ is a divisor on $C$ of degree at least $5$.  This
contradicts B\'ezout's theorem unless $C\subset H$.  But by assumption 
$\text{span}(C)=\PP^4$.  So we conclude that $p\not\in D$.

\ps

Every hyperplane section of $S$ which does not contain $p$ is a
smooth cubic rational curve $D~\subset~ X$.  Thus we have a well-defined 
morphism $I_P~\rightarrow~ \Hdg{3}{0}(X)$.  Let us define $\Pi~\subset~
\HI{2}{F}$ to be the divisor parametrizing planar subschemes $Z~\subset~ 
F$.  Then we can define a morphism 
\begin{equation}
f_P: I_P \rightarrow  \Pi\times \Theta
\end{equation}
as the product of the projection map $I_P\rightarrow \HI{2}{F}$ (which 
factors through $\Pi$ by construction) and the composition of
$I_P\rightarrow \Hdg{3}{0}(X)$ with the Abel-Jacobi map
$\Hdg{3}{0}(X)\rightarrow \Theta$.  

\ps

The claim is that the morphism $f_P$ is injective.  Recall that the
fiber of the Abel-Jacobi map $\Hdg{3}{0}(X)\rightarrow \Theta$
containing some curve $[D]\in \Hdg{3}{0}(X)$ is an open set of the
2-dimensional linear series determined by $D$ on the cubic surface
$X\cap \text{span}(D)$.  The scheme determined by $Z$ will intersect
this cubic surface in a 0-dimensional scheme of length 2.  Such a
scheme imposes 2 linearly independent conditions on divisors in the
linear series $|D|$.  Therefore there is a unique curve $D$ in this
linear system which contains this 0-dimensional scheme of length 2.
Given the curve $D$ and the distinguished point $p$ (which is
determined by $Z$), we can reconstruct the scroll $\Sigma$ as the cone 
over $D$ with vertex $p$.  We can then reconstruct $C$ as the
curve residual to the scheme determined by $Z$ and $D$ in the
intersection $\Sigma\cap X$.  Thus we can uniquely recover $[C]$ from
$f_P([C])$ which shows that $f_P$ is injective.  Therefore
$\dim I_P \leq \dim  \Pi + \dim  \Theta = 3 + 4 = 7$.
This proves the lemma.
\end{proof}

Now define $I_U \subset I$ to be the Zariski dense open subset
parametrizing pairs $([C],[Z])$ with $Z\subset A_C$ a 0-dimensional
scheme of length $2$ such that $\text{span}(C)=\PP^4$ and $Z$ is
non-planar.  If we consider $A_C$ as a subscheme of
$\text{Sym}^2(C)~\cong~\PP^2$, then the length 2 subscheme $Z\subset
\text{Sym}^2(C)$ determines a line in $\text{Sym}^2(C)$, i.e., a linear 
series of degree 2 divisors on $C$.  One consequence of the assumption 
that $Z$ is non-planar is that this linear series has no base-points.
By lemma~\ref{lem-g1240}, there is a unique cubic scroll
$\Sigma\subset \PP^4$ which contains $C$ and such that the linear
series of degree 2 divisors is the linear series of intersections of
$C$ with the lines of the ruling of $\Sigma$.  Let $D$ denote the
directrix of $\Sigma$ and let $F$ denote the divisor class of a line
of the ruling.  Then the hyperplane class on $\Sigma$ is $H~\sim~D+2F$ 
so the intersection $X\cap \Sigma$ has divisor class $3D+6F$.  Now
$C.F=2$ and $C.H = 4$, thus $C~\sim~ 2D+2F$.  On the other hand, the scheme
determined by $Z$ (either $L_1\cup L_2$ if $Z$ is reduced, or the
ribbon $R$ if $Z$ is non-reduced) has divisor class $2F$.  Thus the
residual to $C$ and the subscheme determined by $Z$ is a divisor
$D_2\subset \Sigma$
linearly equivalent to $D+2F$.

\ps

\begin{thm}~\label{thm-irr40} The space $\Hdg{4}{0}(X)$ is
irreducible of dimension $8$. 
\end{thm}

\begin{proof} 
We continue to use the notation introduced in this section.
Because of lemma~\ref{lem-irr40} and lemma~\ref{lem-planar}, it is
equivalent to prove the $I_U$ is irreducible of dimension $8$.  
We stratify $I_U$ according to the type of
the residual curve $D_2$ defined above.  If $D_2$ is a smooth curve,
we say it is the first type.  If $D_2$ is the union of a conic and a
line of the ruling, we say it is the second type.  If $D_3$ is the
union of the directrix $D$ and two lines of the ruling (possibly one
non-reduced line), we say it is the third type.  Define the
corresponding loci in $I_U$ to be $I_1$, $I_2$, and $I_3$.

\ps

\textbf{Third type:}
First we deal with the third type because it is the most involved.  
We will prove that every irreducible component of $I_3$ has dimension
at most $7$.
We 
can associate to each pair $([C],[Z])\in I_3$ the configuration
$([D],[W])$ where $[D]\in F$ is the directrix line and $[W]\in
\HI{4}{F}$ is the length 4 scheme parametrizing the residual to $D$
and $C$ in $\Sigma$.  
Because of lemma~\ref{lem-planar}, we may disregard the subvariety of
$I_3$ such that any length 2 subscheme of $W$ is planar.  Let us
define $I_{3,U}\subset I_3$ to be the open subset such that no length
2 subscheme of $W$ is planar.  Then we are reduced to proving that
every irreducible component of $I_{3,U}$ has dimension at most $7$.

\ps

Notice that $W$ is a subscheme of the divisor
$Z_D\subset F$ which parametrizes lines in $X$ which intersect $D$.
Let $M\subset F\times\HI{4}{F}$ be the closed subset parametrizing
configurations $([D],[W])$ such that $W\subset Z_D$.  Then we have a
morphism $f_3:I_{3,U}~\rightarrow~M$.

\ps

Now we consider the dimension of $M$.  By ~\cite[lemma 10.5]{CG}, for
a general $[D]\in F$, the divisor $Z_D$ is smooth.  Therefore the
fiber of projection onto the first factor $M\rightarrow F$ over the
point $[D]$ is simply the 4-dimensional scheme $\text{Sym}^4(Z_D)$.
It also follows by ~\cite[lemma 10.5]{CG} that when $Z_D$ is singular, 
the singular set consists of a single ordinary node.  The dimension of 
the space of length d subschemes of a nodal curve whose reduced scheme 
is the node is $0$ if $d=1$ and $1$ if $d>0$.  Therefore, even when
$Z_D$ is singular, the dimension of $\HI{4}{Z_D}$ is still only $4$.
So we conclude that $M~\rightarrow~F$ has fiber dimension $4$ so that
$\dim (M) = 6$.

\ps

Now we consider the fiber dimension of $f_3:I_{3,U}~\rightarrow~M$.  We
will prove that every fiber of $f_3$ has dimension at most $1$.
Suppose given a configuration $([D],[W])$.  We want to consider the
quotient $\PP^4/D ~\cong~\PP^2$.  Because no length 2 subscheme of $W$ 
is planar, the image of $W$ in $\PP^4/D$ is a 0-dimensional subscheme
of length 4.  Given any scroll $\Sigma$ which contains $D$ as the
directrix, the image of $\Sigma$ in $\PP^4/D$ is a smooth plane conic.
Given a 0-dimensional subscheme of $\PP^2$ of length 4, there is a
pencil of conics containing this subscheme.  So if $\Sigma$ is a
scroll which contains both $D$ and $W$, then the image of $\Sigma$ in
$\PP^4/D$ must be a conic $B$ in the pencil determined by the image of
$W$.  Moreover, the scroll $\Sigma$ determines an isomorphism
$\phi:D~\rightarrow~B$ which associates to each point in $D$, the
image in $\PP^4/D$ of the line of the ruling through that point.
Notice that given a conic $B$ in $\PP^4/D$ which contains the image of 
$W$, there is at most one isomorphism $\phi:D\rightarrow B$ such that
$\phi$ identifies the projection of $W$ onto $D$ with the projection
of $W$ onto $B$ (because the only automorphism of $\PP^1$ which fixes
a divisor of degree $4$ is the identity map).

\ps

By the last paragraph, associated to a configuration $([D],[W])$ there 
is a pencil of smooth conics $B$ in $\PP^4/D$ which contain the image of $W$, 
and to each $B$ there is (at most) one automorphism $\phi:D~\rightarrow~
B$ which identifies the two projections of $W$.  What extra
information is needed to determine a scroll $\Sigma$ such that
$\Sigma$ contains $D$ and $W$ and such that the image of $\Sigma$ is
$B$?  To answer this question we recall lemma~\ref{lem-rib1} which
says that to determine a scroll $\Sigma\subset \PP^4$ which contains a 
line $D$ as the directrix, it is equivalent to determine the rank 1
subbundle $T\subset T_{\PP^4}|_D$, where $T~\cong~\OO_D(-1)$ is the
bundle of tangent spaces to line of the ruling of $\Sigma$.  Also we
have that $T$ is everywhere distinct (as a rank 1 subspace of
$T_{\PP^4}$) from $T_D$.  So first we consider the image of $T$ in
$N_{D/\PP^4}$.  But of course this is just
$N_{D/\Sigma}~\subset~N_{D/\PP^4}$.  Also $N_{D/\PP^4} =
\OO_D(1)\otimes_\CC N$ where $N$ is the rank 3 vector space whose
associated projective space is canonically $\PP^4/D$.  The map
$\phi:D\rightarrow \PP^4/D$ is equivalent to a subbundle $N'\subset
\OO_D\otimes_\CC N$ where $N'$ is isomorphic to $\OO_D(-2)$,  and the
subbundle $N_{D/\Sigma}\subset \OO_D(1)\otimes_\CC N$ is simply
$N'\otimes \OO(1)$.  So $N_{D/\Sigma}\subset N_{D/\PP^4}$ is uniquely
determined by the map $\phi:D\rightarrow \PP^4/D$.  

\ps

Finally, to determine the scroll $\Sigma$, we have to determine a
subbundle $T\subset T_{\PP^4}|_D$ which projects isomorphically to
$N_{D/\Sigma}$.  The set of such subbundles is equivalent to the set
of global sections of the bundle $\text{Hom}(\OO_D(-1),T_D\oplus
N_{D/\Sigma})~\cong~ \OO_D\oplus \OO_D(3)$ (modulo non-zero scaling)
whose composition with $T_D\oplus N_{D/\Sigma}\rightarrow
N_{D/\Sigma}$ is an isomorphism.  In other words, the set of such
subbundles is simply $\text{Hom}(\OO_D(-1),T_D)~\cong~\OO_D(3)$.  
But this section is determined along the projection of $W$ since $W$
must be a subscheme of the scheme of lines of the ruling of $\Sigma$.
Since a length $4$ subscheme of $\PP^1$ imposes 4 linear conditions on 
$\OO_D(3)$, we see that there is a unique section which restricts to
$W$ in the appropriate way.  So finally we conclude that the scroll
$\Sigma$ is determined by the configuration $([D],[W])$ together with
the conic $B$.  Since $M$ has dimension 6, and since for each
$([D],[W])$ there is at most a 1-dimensional family of possible $B$'s, 
we conclude that $I_{3,U}$ has dimension at most $7$.  

\ps

\textbf{Second Type:}
Next we consider the second type.  We will prove that every
irreducible component of $I_2$ has dimension at most 7.
Let $B$ denote the smooth conic.
Again let $W\subset Z_B\subset F$ denote the length $3$ subscheme
parametrizing the lines which make up the residual to $C\cup B\subset
X$.  By lemma~\ref{lem-planar}, we may suppose that every length
2 subscheme of $W$ is non-planar.  The claim is that the subscheme
parametrized by $W$ spans $\PP^4$.  By way of contradiction, suppose
that it is contained in a hyperplane $H$.  By B\'ezout's theorem, $B$
is also contained in $H$.  But then the intersection of $H$ and the
scroll $\Sigma$ contains the degree 5 curve which is the union of $B$
and the subscheme parametrized by $W$.  This contradicts B\'ezout's
theorem unless $\Sigma~\subset~H$.  But then $C$ is also contained in
$H$, and this contradicts the hypothesis on $C$.  Therefore the scheme 
parametrized by $W$ spans $\PP^4$.  

\ps

Let $M_2\subset \Hdg{2}{0}(X)\times \HI{3}{F}$ be the locally closed subset
parametrizing configurations $([B],[W])$ such that $B$ is smooth, such 
that every length 2 subscheme of $W$ is planar, such that the
subscheme of $X$ parametrized by $W$ spans $\PP^4$, and such that
$W\subset Z_B$, where $Z_B\subset F$ is the locally closed set which
parametrizes lines which intersect $B$ exactly once (there is exactly
one line which intersects $B$ twice).  By the same type of argument at 
in the first case, we conclude that $\dim  M_2 ~=~ \dim 
\Hdg{2}{0}(X) + 3 = 4+3 = 7$.

\ps

There is an obvious morphism $f_2:I_2\rightarrow M_2$, and we are
reduced to showing that this map is injective.  Now given a subscheme
$[W]\in \HI{3}{\text{Grass}(2,5)}$ such that no length 2 subscheme of
$W$ is planar, and such that the scheme parametrized by $W$ (in $\PP^4$) spans
$\PP^4$, then there is precisely one line $L$ whose intersection with
this scheme is of length $3$.
We will only give the
proof when $W$ is reduced -- the non-reduced case is only slightly
more technical.  Suppose that $W=\{[L_1],[L_2],[L_3]\}$.  Then $L_1,
L_2, L_3$ are all disjoint.  If a line $L$ intersects $L_1$ and $L_2$, 
then it lies in the hyperplane $H=\text{span}(L_1,L_2)$.  Since
$\text{span}(L_1,L_2,L_3)=\PP^4$, the line $L_3$ is not contained in
$H$.  Therefore $H\cap L_3$ is a point $p$ which does not lie on $L_1$ 
or $L_2$.  We conclude that the lines $L$ which intersect $L_1,L_2$,
and $L_3$ are exactly the lines $L\subset H$ which intersect $L_1,L_2$ 
and which pass through $p$.  If we consider projection away from $p$,
then the set of such lines corresponds to the intersection points in
$H/p~\cong~\PP^2$ of the images of $L_1$ and $L_2$.  Since these lines 
are skew and don't contain $p$, their images in $H/p$ consist of two
distinct lines,  and two distinct lines in $\PP^2$ intersect in
precisely one point.  

\ps

But given a scroll $\Sigma$, the directrix line $D$ is a line which
intersects $L_1, L_2$ and $L_3$.  Thus we conclude that the directrix
line $D$ is uniquely determined by the configuration $([B],[W])$ (in
fact just by $[W]$).  Moreover, the lines of the ruling induce an
isomorphism $\phi:D\rightarrow B$ which carries each intersection
$L_i\cap D$ to the intersection $L_i\cap B$.  There is a unique
isomorphism $\phi:D\rightarrow B$ with this property (because the only
automorphism of $\PP^1$ which fixes a length 3 divisor is the
identity).  Thus $\phi$ is also determined by the configuration
$([B],[W])$.  From $\phi$ we recover the scroll $\Sigma$.
From the scroll $\Sigma$ we recover $C$ as the residual to $B\cup
L_1\cup L_2\cup L_3$ in $\Sigma\cap X$.  Thus we conclude that $f_2$
is injective, which proves that $I_2$ has dimension at most $7$.

\ps

\textbf{First type:}
Finally we consider $I_1$.  This analysis will also be very important
for describing the Abel-Jacobi map
$\alpha_{4,0}:\Hdg{4}{0}(X)~\rightarrow~ J(X)$.  
Denote the residual curve by $A$ (this is the curve we were calling
$D_2$).
Let $N\subset \Hdg{3}{0}(X)\times F\times F$ denote the locally closed subscheme
parametrizing triples $([A],[L_1],[L_2])$ such that $L_1, L_2$ are
skew lines each of which intersects $A$ transversely in $1$ point and
such that $\text{span}(A,L_1,L_2)~=~\PP^4$.
Of course $N$ fibers over $\Hdg{3}{0}(X)$ and the fiber is an open
subset of $D_A\times D_A$.  Thus we have
\begin{equation}
\dim (N)~=~\dim (\Hdg{3}{0})+2\dim (D_A) ~=~ 6+2~=~8.
\end{equation}
There is an obvious map $I_1~\rightarrow~ N$, and the only nontrivial
condition to verify is that $\text{span}(A,L_1,L_2)~=~\PP^4$, but this follows 
by applying B\'ezout's theorem to $\Sigma$.  

\ps

What are the fibers of
$I_1~\rightarrow~ N$?  Consider the $3$-plane $P~=~\text{span}(L_1,L_2)$.
This intersects $A$ in a degree $3$ divisor.  Two of the points of
this divisor are the points of intersection of $A$ and $L_1$, $L_2$.
The third point $p$ lies on neither $L_1$ nor $L_2$ since $A.L_i$ is a
degree $1$ divisor.  Now there is a unique line $M$ which contains $p$ 
and which intersects both $L_1$ and $L_2$: if we project $P$ away from 
$p$, then the line $M$ corresponds to the unique point of intersection 
of the images of $L_1$ and $L_2$ in $\PP^2$.  Now suppose that
$\Sigma$ is a scroll which contains $L_1$ and $L_2$ and $A$.  Let $D$
denote the directrix.  Since $D$ intersects $L_1$ and $L_2$, it must
lie in $P$.  If $D$ does not contain $p$, then there is a line of the
ruling $F$ of $\Sigma$ which passes through $p$.  But then $L_1\cup
L_2\cup D\cup F$ is a divisor of degree $4$ in the hyperplane section
$P\cap \Sigma$.  This contradicts B\'ezout's theorem.  What we
conclude is that $D$ must equal $M$.  Moreover the isomorphism
$\phi:A~\rightarrow~ M$ corresponding to projection of
$\Sigma~\rightarrow~ M$ must be the unique isomorphism such that
$\phi(p)~=~p$ and such that $\phi(A\cap L_i) ~=~ \phi(M\cap L_i)$.  Of
course the scroll $\Sigma$ is determined by $M$ and the isomorphism
$\phi$.  Thus there is a unique scroll $\Sigma$ which contains $A\cup
L_1\cup L_2$, i.e., $I_1~\rightarrow~ N$ maps $1$-to-$1$ to its image.
Since there is some irreducible component of $I$ of dimension $8$, in
fact we must have that $I_1~\rightarrow~ N$ is dominant,
i.e., $I_1~\rightarrow~ N$ is an open immersion.  Finally notice that $N$ 
fibers over the irreducible space $\Hdg{3}{0}$ and the general fiber
$D_A\times D_A$ is irreducible.  Thus $N$ is irreducible.  So $I_1$ is 
irreducible of dimension $8$.

\ps

Since $I\rightarrow \Hdg{4}{0}(X)$ is surjective, every 
component of $\Hdg{4}{0}(X)$
is dominated by a component of $I$, which must be at least eight dimensional,
since every component of $\Hdg{4}{0}$ is.  The only component of $I$ with this
property is $\overline{I_1}$, which is precisely eight dimensional, so we 
conclude that $\Hdg{4}{0}$ is irreducible of dimension eight.
\end{proof}

\textbf{Remark:}
Let $I_1$ be as in the proof above and let $J_1$ be the quotient of
$I_1$ by the involution $([C],[L_1],[L_2])~\mapsto~ ([C],[L_2],[L_1])$.
Notice that $J_1~\rightarrow~ \Hdg{4}{0}(X)$ is still dominant.

\section{Quintic Elliptic Curves}\label{sec-irr51}

In section~\ref{sec-irr40} we proved that $\Hdg{4}{0}(X)$ is
irreducible by residuating the union of a quartic curve and a pair of
2-secant lines in the intersection of $X$ with a suitable cubic scroll 
$\Sigma$.  In this section we will prove that $\Hdg{5}{1}(X)$ by
residuating a quintic genus 1 curve in the intersection of $X$ with a
suitable cubic scroll.  The idea of the proof is very similar to the
proof of theorem~\ref{thm-irr40}.  As in that proof, there are several 
degenerate behaviors which we need to rule out as generic.

\ps

\begin{thm}~\label{thm-irr51} 
The space $\Hdg{5}{1}(X)$ is irreducible of dimension $10$.
\end{thm}

\begin{proof}
By lemma~\ref{lem-equiv1} the irreducibility of $\Hdg{5}{0}$ is equivalent 
to showing that
$\widetilde{H}^{5,1}$ is irreducible, where
$\widetilde{H}^{5,1}$ is the parameter space for pairs
$(f:\Sigma~\rightarrow~ \PP^4, i:E~\rightarrow~ \Sigma)$ such that $f\circ
i:E~\rightarrow~ \PP^4$ is an embedding of $E$ as a quintic elliptic
curve.  Indeed, we have seen that the fiber of projection
$\widetilde{H}^{5,1} ~\rightarrow~ \Hdg{5}{1}(X)$ over a point $[E]$
is simply the set of $g^1_2s$ 
on $E$, i.e., $\text{Pic}^2(E) \cong E$.  Since the fibers are
irreducible of constant fiber dimension 1, we conclude that $\Hdg{5}{1}(X)$
is irreducible of dimension $10$ iff $\widetilde{H}^{5,1}$ is
irreducible of dimension $11$.  

\ps

On the other hand each pair $(f:\Sigma~\rightarrow~
\PP^4,i_E:E~\rightarrow~\Sigma)$ is equivalent to a pair
$(f:\Sigma~\rightarrow~ \PP^4, i_C:C~\rightarrow~ \Sigma)$ where 
$C$ is the residual quartic curve, $[C]=D+3F$.  
We decompose $\widetilde{H}^{5,1}$ into
a union of locally closed subsets
$\widetilde{H}_1$, $\widetilde{H}_2$, $\widetilde{H}_3$, $\widetilde{H}_4$
parametrizing the 
set where $i_C:C~\rightarrow~ \Sigma$ is in the first, second, third or
fourth case (we say that $i_C:C~\rightarrow~ \Sigma$ is in the $i$-th case
if $C'.D ~=~ i-2$ where $C'$ is the unique irreducible component of $C$
which projects isomorphically to $\PP^1$ under $\pi$).  We will show that
for $i\neq 1$, $\widetilde{H}_i$ has dimension $\leq 10$, and we will
show that $\widetilde{H}_1$ is irreducible of dimension $11$.  

\ps

\textbf{First Case:}
Now $\widetilde{H}_1$ parameterizes pairs ($f:\Sigma\rightarrow\PP^4$,
$i_C:C\rightarrow\Sigma$) where $\Sigma$ is in the first case, and
$f(i_C(C))\subset X$.  There is a projection $\widetilde{H}_1 \rightarrow
\Hdg{4}{0}$ which assigns to $(f,i_C)$ the curve $C\subset X$, the
embedding being given by 
$f\circ i_C$.  We have seen in lemma \ref{lem-scr40} that the fibre 
of this projection over a particular curve $C\subset \PP^4$ consists of
the data of a line $L$ in $\PP^4$ intersecting $C$ in a subscheme $Z$ of
length $2$, and an isomorphism $\phi$ between $C$ and $L$ which is the
identity map on $Z$.  The length $2$ subscheme $Z$ uniquely determines
$L$, and given a fixed $Z$, there is a $\CC^{*}$ worth of choices of 
such isomorphisms $\phi$.  Therefore each fibre of $\widetilde{H}_1
\rightarrow\Hdg{4}{0}$ is itself a $\CC^{*}$ bundle over the space
$\Sym^2(C)=\PP^2$ parameterizing the $Z$'s.  We see that $\widetilde{H}_1
\rightarrow \Hdg{4}{0}$ has irreducible fibres of dimension $3$.  
By theorem~\ref{thm-irr40},  $\Hdg{4}{0}(X)$ is irreducible of dimension $8$,
and therefore $\widetilde{H_1}$ is irreducible of dimension $11$.

\ps

\textbf{Second Case:}
Now $\widetilde{H}_2$ parametrizes pairs $(f:\Sigma~\rightarrow~
\PP^4$, $i_C:C~\rightarrow~ \Sigma)$ where $C=C'\cup F$ the union of a
smooth rational cubic curve $C'$ and a line of ruling $F$, and
$f(i_C(C))\subset X$.  Consider the morphism $\widetilde{H}_2~\rightarrow~
H^{3,0}$ which associates to $(f,i_C)$ the curve $C'\subset X$, 
Recall that $\dim (H^{3,0}) ~=~ 6$.  We analyze the fibre of this map
by looking for the data necessary to reconstruct $\Sigma$.
The irreducible component $F\subset C$ is mapped to a line in $X$ which
intersects $C'$ in a single point.  The directrix $D$ of $\Sigma$ is mapped
to a line in $\PP^4$ which intersects $F$, and also intersects $C'$ in a 
single point. Given a fixed $C'\subset X$, there is a one parameter family
of lines in $X$ to serve as an $F$.  Given a fixed $F$, we recover the directrix
as follows: pick any point $p$ on $C'$,  then there is a $\P^1$ of lines $D$
passing through $p$ and intersecting $F$ (in case $p\in F\cap C'$, the limiting
condition is that $D$ lie in the $\PP^2$ spanned by $F$ and the tangent line
to $C'$ at $p$). Finally, we need to specify the isomorphism $\phi:C'\rightarrow
D$ induced by the lines of ruling.  Since this must be the identity on $p$ and
on $F\cap C'$, this is parameterized by $\CC^{*}$.  As in the other lemmas
on reconstructing cubic scrolls in section \ref{sec-cuscr}, this data is
sufficient to specify $\Sigma$.  Altogether we see that the dimension
of $\widetilde{H}_2$ is the sum of $6$ for $\dim(\Hdg{3}{0})$, $1$ for
the choice of the line $F$, $1$ for the choice of point $p\in C'$, 
$1$ for the choice of $D$ going through $p$ and intersecting $F$, and $1$ 
for the $\CC^{*}$ of isomorphisms between $C'$ and $D$ satisfying our
conditions, i.e., $\dim(\widetilde{H}_2)=10$.

\ps

\textbf{Third Case:}
This time the curve $C'$ is a smooth conic, and $C$ consists of $C'$
and two lines $F_1$, $F_2$ of ruling (possibly a double line).  The inclusion
$i_C$ takes 
the lines of ruling to two lines (or possibly a nonplanar ribbon) in $X$ 
which intersect $C'$.  The directrix $D$ of $\Sigma$ maps to a line in $\PP^4$ which
intersects $C'$ once and the union of the lines in a subscheme of length two.  
We have a projection $\widetilde{H}_3~\rightarrow~ H^{2,0}$, given by forgetting
all of the data except the conic $C'$.  Reversing this procedure, if we start
with a smooth conic $C'\subset X$, the choices of two lines $F_1$, $F_2$, in 
$X$ meeting 
$C'$ form a two dimensional family.  Given the two lines, the directrix $D$
must meet each of them, and so is also parameterized by a two dimensional 
family, namely the choices of the intersection points on the two lines.
Finally, given this data, we have to specify the isomorphism $\phi:C'\rightarrow
D$ corresponding to the lines of ruling.  This isomorphism must take 
$F_i\cap C'$ to $F_i\cap D$ for $i=1,2$, and so we see that there is a 
$\CC^{*}$ of choices.
  Altogether the dimension of $\widetilde{H}_3$
is the sum of $4=\dim (H^{2,0})$, $2$ for the union of two lines
intersecting $C'$, $2$ for the 2-parameter family of possibilities for the
directrix $D$, and $1$ for the $\CC^*$ of isomorphisms $\pi:C'~\rightarrow~ D$,
i.e., $\dim (\widetilde{H}_3) = 9$.  

\ps

\textbf{Fourth Case:}
Finally we consider the fourth case.  This time $C'$ is the directrix 
of $\Sigma$, and $C'~\subset~ X$
is a line in $X$.  The Fano scheme of lines
in $X$ has dimension $2$.  The remaining components of $C$ are mapped
to a union of three lines intersecting $C'$ (or some degeneration
thereof).  For fixed $C'$, the dimension of such triples of lines is
$3$.   By lemma \ref{lem-rib1}, in order to construct a scroll $\Sigma$
containing $C'$ as the directrix, we need to provide a sub-line-bundle
$T\subset T_{\PP^4}|{C'}$, with $T\cong\OO_{C'}(-1)$.  The set of such
bundles is a $\PP^{12}$, since $\hom(\OO_{C'}(-1),T_{\PP^4}|_{C'})=13$.
In order for the scroll to contain the three lines touching $C'$, this
sub-bundle must agree with the direction of each line at the point of contact
with $C'$.  For each line, this is a three dimensional linear condition. 
Therefore the space of scrolls containing $C'$ as the directrix, as well as
the three lines as lines of ruling is a $\PP^3$.
Thus altogether $\widetilde{H}_4$ has dimension $(2+3+3)-1~=~7$.

\ps

\end{proof} 

\textbf{Remark}  Of course the proof shows more than just that
$\Hdg{5}{1}(X)$ is irreducible.  We see that for a general quintic elliptic
$E\subset X$ and a general cubic scroll containing $E$, the residual
curve is a smooth quartic rational curve.  

\section{Quintic Rational Curves}

In this section we will prove that the space $\Hdg{5}{0}(X)$ is
irreducible.
\ps

\begin{lem}~\label{lem-rncsec}  Let $C\subset\PP^n$ be a rational
normal curve and let $P\subset \PP^n$ be a linear $r$-plane.
If $(r+2)k\geq (r+1)(n+1)$, then $P$ is contained in a
$k$-secant $(k-1)$-plane of $C$, i.e., there
exists a divisor $D~=~q_1+\dots q_k$ on $C$ such that
$P\subset \text{span}(D)$.
\end{lem}

\begin{proof}
We identify $C$ with $\PP^1$ so that $\OO_C(1)$ is a degree 1 line
bundle, and $\OO_{\PP^n}(1)|_C$ is a degree $n$ line bundle.  Up to a
choice of basis of $\PP^n$, we can identify the inclusion
$C\hookrightarrow \PP^n$ with the morphism associated to the complete
linear series $|\OO_C(n)|$.  

\ps

Let $\PP^k$ be identified with the complete linear series
$|\OO_C(k)|$.  Then on $\PP^k$ we have the tautological injection of
vector bundles $\OO_{\PP^k}(-1)~\rightarrow~ H^0(C,\OO_C(k))\otimes_\CC
\OO_{\PP^k}$.  If we take the tensor product of this map with
$H^0(C,\OO_C(n-k))$ and then use the product map
\begin{equation}
\begin{CD}
H^0(C,\OO_C(n-k))\otimes_\CC H^0(C,\OO_C(k)) @>>>
H^0(C,\OO_C(n)),
\end{CD}
\end{equation}
we have the composite map
\begin{equation}\begin{CD}
H^0(C,\OO_C(n-k))\otimes_\CC \OO_{\PP^k}(-1) @>>>
H^0(C,\OO_C(n))\otimes\OO_{\PP^k}.
\end{CD}\end{equation}
If we think of $\PP^k$ as the parameter space for degree $k$ divisors
$D~=~q_1+\dots + q_k$ on $C$, i.e., as $\Sym^k(C)$, then the fiber 
of this map of vector bundles at a point $[D]$ is just
$H^0(C,\OO_C(n)(-D))~\rightarrow~ H^0(C,\OO_C(n))$.  

\ps

Under the identification $H^0(C,\OO_C(n))~=~H^0(\PP^n,\OO_{\PP^1}(1))$,
we have a restriction map $H^0(C,\OO_C(n))~\rightarrow~
H^0(P,\OO_{\PP^n}(1)|_P)$.  Thus we have an induced map of vector bundles
on $\PP^k$ obtained as the composite map
\begin{equation}\begin{CD}
H^0(C,\OO_C(n-k))\!\otimes_\CC\! \OO_{\PP^k}(-1) @>>>
H^0(C,\OO_C(n))\!\otimes_\CC\! \OO_{\PP^k} @>>>
H^0(P,\OO_{\PP^n}(1)|_P)\!\otimes_\CC\! \OO_{\PP^k}.
\end{CD}\end{equation}
Suppose the fiber of this map is the zero map at a point $[D]$.  Then
every linear polynomial of $\PP^n$ which vanishes on $D$ also vanishes
on $P$.  Since $\text{span}(D)\subset \PP^n$ is cut out by the linear
polynomials which vanish on $D$, we conclude that the ideal of
$\text{span}(D)$ is contained in the ideal of $P$, i.e., $P\subset
\text{span}(D)$.  So we are reduced to showing that some fiber of this
map is zero, i.e., this map of vector bundles has nonempty zero locus.

\ps

We may think of the map above as a global section of the bundle
\begin{equation}
\text{Hom}_\CC(H^0(C,\OO_C(n-k)),H^0(P,\OO_{\PP^k}(1)|_P))\otimes_\CC
\OO_{\PP^k}(1). 
\end{equation}
The rank of this vector bundle is 
\begin{equation}
\dim H^0(C,\OO_C(n-k))\times \dim H^0(P,\OO_{\PP^k}(1)|_P)
~=~ (n+1-k)(r+1).
\end{equation}
Thus the map is a global section of $\OO_{\PP^k}(1)^{(n+1-k)(r+1)}$.
The zero locus is just defined by the vanishing of $(n+1-k)(r+1)$
linear polynomials.  So long as $(n+1-k)(r+1)\leq k$, these linear
polynomials always have a solution.  Thus if $(r+2)k\geq (n+1)(r+1)$,
then the zero locus is nonempty.
\end{proof}

\ps

\begin{Remark}~\label{rem-rncsec} Notice that the proof also shows
that the set of 
$k$-secant $(k-1)$-planes which contain $P$ is a linear subspace of
$\PP^k$.  In particular, when this set is finite, there is a unique
solution.  
\end{Remark}

\ps

\begin{cor}~\label{cor-3sec50}  If $C\subset \PP^4$ is a smooth,
nondegenerate quintic rational curve, then $C$ 
has a unique $3$-secant line $L\subset \PP^4$, and $L$ is not a
$4$-secant line.  If $C\subset \PP^3$ is a smooth quintic rational
curve, the $C$ has a 1-parameter family of $3$-secant lines $L\subset
\PP^3$.  If every $3$-secant line to $C$ is a $4$-secant line, then
$C$ lies on a smooth quadric surface as a divisor of type $(1,4)$.
\end{cor}

\begin{proof}
First consider the case where $C$ is nondegenerate.  Then we can think
of $C\subset \PP^4$ as the projection of a rational normal curve
$C'\subset \PP^5$ from a point $p$ not on $C'$.  
By lemma~\ref{lem-rncsec}, we see that there is a $3$-secant $2$-plane
$\text{span}(D)$ which contains $p$.  The projection of $P$ is a
$3$-secant line $L$ to $C$.  On the other hand, suppose that $C$ has a
4-secant line $L$.  The preimage of $L$ is a $4$-secant $2$-plane to $C'$.
But since any $4$ points on $C$ are linearly independent (or more
generally any degree $4$ divisor on $C$ imposes $4$ conditions on
linear forms), we see that $C'$ does not have a $4$-secant $2$-plane.
Thus $C$ has a $3$-secant line, but does not have a $4$-secant line.

\ps

Suppose that $C$ has two distinct $3$-secant lines $L,M$.  Consider
$H~=~\text{span}(L,M)$.  If this is a hyperplane in $\PP^4$, then $H\cap
C$ has degree $6$.  This contradicts B\'ezout's theorem unless
$C\subset H$, i.e., $C$ is degenerate.  If $H$ is a $2$-plane, choose
any point $p\in C$ not contained in $H$ and let
$H'~=~\text{span}(H,p)$.  Then $H'$ is a hyperplane, and again $H'\cap C
\supset \{p\}\cup (H\cap C)$ has degree at least $6$.  Again by
B\'ezout's theorem we conclude that $C\subset H'$ so that $C$ is
degenerate.

\ps

Suppose now that $C$ is degenerate.  Since $C$ is smooth,
$\text{span}(C)$ is a hyperplane in $\PP^4$.  Thus we may think of $C$
as the projection of a rational normal curve $C'\subset \PP^5$ from a
line $N\subset \PP^5$. 
Now the $3$-secant lines to $C$ correspond to $3$-secant $3$-planes to
$C$ in $\PP^5$ which contain $N$.  It is a bit simpler to think of
this as the set of $3$-secant $2$-planes which intersect the line
$N$.  Since there is one such $2$-plane for each point of $N$, we see
that $C$ has a pencil of $3$-secant lines.  Suppose
moreover that each of these $3$-secant lines is actually a $4$-secant
line.  If two of these lines, $L,M$ intersect nontrivially, then
$P~=~\text{span}(L,M)$ is a $2$-plane and $P\cap C$ has degree at least
$7$.  This contradicts B\'ezout's theorem unless $C\subset P$, which
itself contradicts that $C$ is smooth.  Thus all of the $4$-secant
lines are skew.  Now let $S$ be the surface swept out by the
$4$-secant lines.  Then $S$ contains $C$.  Choose any 2-secant line
$M$ to $C$.  For each $4$-secant line $L$ to $C$ which intersects $M$, 
consider the 2-plane $\text{span}(L,M)$.  If $L$ does not pass through 
one of the 2 points of intersection of $M\cap C$ then
$\text{span}(L,M)$ intersects $C$ in at least $6$ points, which
contradicts B\'ezout's theorem.  Therefore the only lines $L$ which
intersect $M$ are the lines through the 2 points of intersection of
$M\cap C$.  Thus $S$ intersects $M$ in exactly 2 points, i.e., $S$ is a 
quadric surface.  Since $S$ contains a 1-parameter family of skew
lines, we conclude that $S$ is a smooth quadric surface.  Finally,
every smooth quintic rational curve on a smooth quadric surface has
divisor class $(1,4)$ (with respect to some ordering of the two
rulings).
\end{proof}

Now suppose that $[C]\in \Hdg{5}{0}(X)$.  If $C$ has a 1-parameter
family of $4$-secant lines, then we see by corollary~\ref{cor-3sec50}
that $C$ is a divisor of type $(1,4)$ on a smooth quadric $Q$.  But
$Q\cap X$ is a divisor of type $(3,3)$ on $Q$, it cannot contain a
divisor of type $(1,4)$ as an irreducible component.  This
contradiction shows there are no such curves.

\ps

Define $I=I_{5,0}\subset \Hdg{5}{0}(X)\times \mathbb{G}(1,4)$ to be the locally
closed subvariety parametrizing pairs $(C,L)$ where $L$ is a
$3$-secant line to $C$ which is not a $4$-secant line.  Given such a
pair, let $Z ~=~ L\cap C$.  This is a degree 3 divisor on both $L$ and
$C$,  so there is a unique isomorphism $\phi:L~\rightarrow~ C$ such
that $\phi(Z)~=~Z$.  By lemma~\ref{lem-scr50} associated to the
data $C$, $L$, and $\phi$ there is a unique triple $(f,i,j)$ with
$f:\Sigma\rightarrow \PP^4$ a cubic scroll, $i:C\rightarrow \Sigma$ and
$j:L\rightarrow \Sigma$ factorizations of $C\rightarrow \PP^4$, and 
$L\rightarrow \PP^4$, and such that $L$ is the directrix of $\Sigma$.

\ps

Conversely, given a cubic scroll $f:\Sigma~\rightarrow~ \PP^4$ and a
factorization $i:C~\rightarrow~ \Sigma$ of the inclusion with $i(C)\sim
D+4F$, we see that $f(D)$ is a $3$-secant line which is not a
$4$-secant line.  Therefore $I$ also parametrizes triples
$(C,f:\Sigma~\rightarrow~ \PP^4,\phi)$.

\ps

Now for each cubic scroll $f:\Sigma~\rightarrow~ \PP^4$ and
$j:C~\rightarrow~ \Sigma$ as above, the residual $C_2$ to $j(C)$
in $f^{-1}(X)$ is a divisor of type $2D+2F$.  We know from
subsection~\ref{subsec-scr40} that such a divisor is a quartic curve
of arithmetic genus $0$, e.g. a quartic rational curve.

\ps

\begin{thm}~\label{thm-irr50}  $\Hdg{5}{0}(X)$ is irreducible of
dimension 10.  For a 
general $[C]\in 
\Hdg{5}{0}(X)$, if $f:\Sigma~\rightarrow~ \PP^4$ is the unique cubic
scroll containing $C$, the residual curve $C_2$ to $C\subset f^{-1}(X)$ is a 
smooth quartic rational curve. 
\end{thm}

\begin{proof}
Decompose $I$ depending
on the \emph{type} of $C_2$.  We say $C_2$ is the first type if
it is a smooth quartic rational curve.  We say $C_2$ is the
second type if $C_2$ is a union of two smooth conics $A\cup B$.  We
say that 
$C_2$ is the third type if it is a union of the directrix and a
twisted cubic $D\cup A$.  We say that $C_2$ is the fourth type if it
is the 
union of a conic, the directrix, and a line of the ruling  $A\cup D\cup
F$.  We say
that $C_2$ is the fifth type if $C_2$ is the union of the double of
the directrix and two lines of the ruling $2D\cup F_1\cup F_2$.  Finally,
we say that $C_2$ is 
of the sixth type if $C_2$ is the double of a conic.  We will label
the corresponding locally closed subsets of $I$ by $I_1,\dots, I_6$.

\ps

First we show that for each $i> 1$, $\dim I_i \leq 9$.

\ps

\textbf{Second type:}
Suppose
that $C_2$ is the second type.  The scroll $f:\Sigma~\rightarrow~
\PP^4$ is determined by giving the union of the two conics $A\cup
B$ meeting at a point $p$, and by giving the isomorphism
$\phi:A~\rightarrow~ B, \phi(p)~=~p$ induced by the lines of the ruling of
$\Sigma$.  Thus we see that $I_2$ fibers over the Hilbert scheme of
intersecting conics with fibers of dimension $2$: the set of
isomorphisms is a principal homogeneous space for the $2$-dimensional
subgroup of automorphisms in $PGL(2)$ which fix a point of $\PP^1$.
To specify a conic
in $X$ it is equivalent to specify a line in $X$ and a $2$-plane
containing this line (the conic is the residual of the line).  Thus to
specify two conics intersecting in a point $p\in X$, it is equivalent
to specify a pair of lines $L,M$ and then let the $2$-planes be
$\text{span}(L,p)$ and $\text{span}(M,p)$.  So we see that the Hilbert
scheme of intersecting conics is birational to $X\times
\text{Sym}^2(F)$, and so has dimension $3+2+2~=~7$.  So $I_2$ has
dimension $2+7~=~9$.  

\ps

\textbf{Third type:}
Suppose that $C_2$ is the third type.  To specify the scroll it is
equivalent to specify the twisted cubic $A$, the directrix $D$ which
intersects $A$ in a point $p$, and an isomorphism $\phi:A~\rightarrow~
D$ such that $\phi(p)~=~p$.  Thus $I_3$ fibers over the Hilbert scheme
of unions $A\cup D$ with fibers which are $2$-dimensional.  We have
seen that $\Hdg{3}{0}$ has dimension $6$,  and that the set of lines
intersecting a twisted cubic $A$ has dimension $1$.  Thus the Hilbert
scheme of unions $A\cup B$ has dimension $7$.  So $I_3$ has dimension
$2+7~=~9$.

\ps

\textbf{Fourth type:}
Suppose that $C_2$ is the fourth type.  To specify the scroll it is
equivalent to specify the directrix line $D$, the conic $A$, a line of
ruling $F$ intersecting both $D$ and $A$ (in distinct
points), and an isomorphism $\phi:D~\rightarrow~ A$ such that
$\phi(F\cap D) ~=~ F\cap A$.  Thus $I_4$ fibers over the Hilbert scheme
of curves $A\cup D\cup F$ with fibers which are $2$-dimensional.  To
specify $A\cup D\cup F$, it is equivalent to specify $D\cup F$, a
point $p\in F$ and the residual line $L$ to $A$.  The dimension of the
space of intersecting lines is $3$.  The dimension of choices for $p$
is $1$,  and the dimension of choices for $L$ is $2$.  Thus the
dimension of the space of curves $A\cup D\cup F$ is $3+1+2~=~6$.  So
$I_4$ has dimension $8$.

\ps

\textbf{Fifth type:}
Suppose that $C_2$ is the fifth type.  By lemma~\ref{lem-rib1}, we
know that the scroll $\Sigma$ is determined by the vertical tangent
bundle $T\subset T_{\PP^4}|_D$.  The condition that the intersection
of $\Sigma$ contain the double of $D$ is exactly that the normal
bundle $N_{D/\Sigma}\subset N_{D/\PP^4}$ is contained in $N_{D/X}$.
But this normal bundle is simply the image of $T$ in the quotient
$N_{D/\PP^4}$ of $T_{\PP^4}|_D$.  Thus the scrolls $\Sigma$ such that
$\Sigma\cap X$ contains $2D$ are the same as sub-line-bundles $T\subset
T_X|_D$ of degree $-1$.  In both of the cases $T_X|_D \cong \OO_D(2)\oplus
\OO_D\oplus \OO_D$ and $T_X|_D\cong \OO_D(2)\oplus
\OO_D(1)\oplus\OO_D(-1)$ we have that $H^0(C,T_X|_D(1))$ is an
$8$-dimensional vector space.  
Moreover, each of the two lines $F_1$, $F_2$ of the ruling contained
in $X$ imposes two
linear conditions on the
sections.  Since the set of scrolls is the projective space associated to
the possible sections, we see that there is at most a $3$-dimensional family 
of scrolls which contain the double of $D$ and $F_1$, $F_2$.  Therefore
the dimension of the
space of pairs $([D],\{F_1,F_2\})$ is just the sum of $2$ for the line in $X$,
and $1$ each for the $F_i$'s.   Altogether,
we see that $\dim (I_5)~\leq~2+1+1+3~=~7$.

\ps

\textbf{Sixth type:}
Finally we consider the sixth type.  By lemma~\ref{lem-rib2}, to
specify a scroll containing a conic $A\subset X$ is the same as giving
a sub-line bundle $T\subset T_{\PP^4}|_A$ of degree $1$.  As in the
last case, the condition that $\Sigma\cap X$ contain $2A$ is exactly
that $T\subset T_X|_A$.  The two possibilities for $T_X|_A$ are
$\OO_A(2)\oplus \OO_A(1)\oplus \OO_A(1)$ and
$\OO_A(2)\oplus\OO_A(2)\oplus\OO_A$.  In both cases we see that
$H^0(A,T_X|_A(-1))$ is a vector space of dimension $4$.  Thus there is
a $\PP^3$ of scrolls $\Sigma$ such that $\Sigma\cap X$ contains $2A$.
Since $\dim (\Hdg{2}{0}(X))~=~4$, we conclude that $\dim (I_6)
~=~ 4+3 ~=~ 7$.

\ps

In each case the dimension is at most $9$.  
By proposition~\ref{prop-def1} every irreducible component of
$\Hdg{5}{0}$ has dimension at least $10$.  By
corollary~\ref{cor-3sec50} we know $I\rightarrow \Hdg{5}{0}$ is
surjective.  So we conclude that the image of $I_1\rightarrow
\Hdg{5}{0}$ is Zariski dense and has dimension at least $10$.

\ps

Fixing a quartic rational curve
$C_2\subset X$, by lemma~\ref{lem-g1240} the set of cubic scrolls
$f:\Sigma~\rightarrow~ \PP^4$ containing 
$C_2$ is equivalent to the set of (basepoint free) $g^1_2$'s on $C_2$.
The set of 
$g^1_2$'s on $C_2$ is simply $\text{Sym}^2(C_2)~\cong~\PP^2$.
We see that $I_1$ fibers over $\Hdg{4}{0}$ as a
$\PP^2$-fibration.  By theorem~\ref{thm-irr40}, $\Hdg{4}{0}$ is
irreducible of dimension $8$.  Thus $I_1$ is irreducible of dimension $10$.
So the image of $I_1\rightarrow \Hdg{5}{0}$ is irreducible of
dimension at most $10$.  On the other hand we know the image has
dimension at least $10$.  So $\Hdg{5}{0}$ is irreducible of dimension $10$.
\end{proof}

\section{Quintic Curves of Genus $2$}

By B\'ezout's theorem, $X$ cannot contain a plane curve of degree
$d>3$.  Thus the next case after quintic elliptic curves is quintic
curves of genus $2$.  

\ps

Suppose $C\subset X$ is a quintic curve of genus
$2$.  Let $H$ denote the hyperplane class on $C$.  Since
$\text{deg}(H)~=~5 > 2 ~=~ \text{deg}(K_C)$, we conclude that
$H^1(C,\OO_C(H))~=~0$.  Thus by Riemann-Roch we have
\begin{equation}
\dim H^0(C,\OO_C(H)) ~=~ \text{deg}(H)+1-g~=~5+1-2~=~4.
\end{equation}
Thus the complete linear system $|H|$ is a $\PP^3$, i.e., $C$ is
contained in a $\PP^3$ inside $\PP^4$.  Moreover, by Riemann-Roch we
also have that 
\begin{equation}
H^0(C,\OO_C(2H))~=~10+1-2~=~9 < 10 ~=~ H^0(\PP^3,\OO_{\PP^3}(2)).
\end{equation}
Therefore $C$ is contained in a quadric surface $C\subset S$.  Now
$S\cap X$ is a Cartier divisor of degree $6$ on $S$.  Since $C$ is
degree 5, the residual of $C\subset S\cap X$ is a divisor of degree 1,
i.e., a line.  Therefore every quintic genus $2$ curve is residual to a
line $L\subset X$ in a quadric surface.  

\ps

Let $\PP(Q^\vee) ~\rightarrow~ F$ denote the $\PP^2$-bundle over $F$
parametrizing pairs $([L],[H])$ where $L\subset H\subset \PP^4$ is a
line contained in a hyperplane contained in $\PP^4$ such that
$L\subset X$.  Let $U~\rightarrow~ \PP(Q^\vee)$ denote the $\PP^6$-bundle
parametrizing triples $([L],[H],[S])$ where $S\subset H$ is a quadric
surface containing $L$.  The universal quadric surface
$\widetilde{S}\subset U\times 
\PP^4$ is a Cartier divisor 
inside the pullback of the universal hypersurface $\widetilde{H}\subset
U\times \PP^4$.  Since $U$ is smooth so is $\widetilde{H}$, therefore
$\widetilde{S}$ is a local complete intersection.  Next, $U\times X\subset
U\times \PP^4$ is a Cartier divisor.  Since $\widetilde{S}$ and $U\times
X$ have no irreducible component in common, we see that $D:~=~\widetilde{S}\cap
U\times X \subset \widetilde{S}$ is a Cartier divisor locally cut out by a
regular element, so $D$ is also a local complete intersection.  In
particular, $D$ is Gorenstein.  

\ps

Let $D_1\subset D$ denote the pullback from $F$ of the universal line
in $X$.  Then $D_1~\rightarrow~ U$ is smooth, therefore $D_1$ is
smooth.  In particular $D_1$ is Cohen-Macaulay.  Therefore by
corollary~\ref{cor-reform}, the residual $D_2$ to $D_1$ in $D$ is a
flat family of 
Cohen-Macaulay schemes.  By specializing to a point $([L],[H],[S])$
with $S$ smooth, we see that the general fiber of $D_2$ is a smooth
quintic genus 2 curve in $X$.  Thus there is an induced map
$f:U~\rightarrow~  \Hi{5t-1}(X)$.  We have seen that this
map is a bijection over $\Hdg{5}{2}(X)$.  Therefore the preimage
$f^{-1}(\Hdg{5}{2}(X))$ is precisely the normalization
$\widetilde{\Hdg{5}{2}}(X)$.  Since $U$ is irreducible, we also 
conclude that $f(U)~=~\overline{\Hdg{5}{2}}(X)$.  Thus we have the
following result.

\ps

\begin{thm}  The normalization $\widetilde{\Hdg{5}{2}}(X)$ of
$\Hdg{5}{2}(X)$ is a smooth, connected variety of dimension $10$.  
\end{thm}

\bibliography{my}
\bibliographystyle{abbrv}

\end{document}